\documentclass[aop,seceqn,MSNbibl,citesort,dvips]{arximspdf}
\usepackage{mathbh}
\usepackage{graphicx}
\usepackage{breakurl}

\doi{10.1214/10-AOP539}
\volume{38}
\issue{6}
\pubyear{2010}
\firstpage{2379}
\lastpage{2417}

\makeatletter

\newtheorem{theorem}{Theorem}[section]
\newtheorem{lemma}[theorem]{Lemma}
\newtheorem{prop}[theorem]{Proposition}
\newtheorem{cor}[theorem]{Corollary}

\newproclaim{defin}[theorem]{Definition}
\newproclaim{remarks}[theorem]{Remarks}

\newproclaim{Remark}{Remark}

\newcommand{\disjoint}{\cap\hspace*{-9pt} / \hspace*{3.5pt}}

\makeatother

\begin{document}
\begin{frontmatter}

\title{Exponential tail bounds for loop-erased random walk in two dimensions}
\runtitle{Exp tail bounds for LERW in 2D}

\begin{aug}
\author[A]{\fnms{Martin T.} \snm{Barlow}\thanksref{t1,t2}\ead[label=e1]{barlow@math.ubc.ca}} and
\author[A]{\fnms{Robert} \snm{Masson}\corref{}\thanksref{t1}\ead[label=e2]{rmasson@math.ubc.ca}}
\runauthor{M. T. Barlow and R. Masson}
\affiliation{University of British Columbia}
\address[A]{Department of Mathematics\\
University of British Columbia \\
Vancouver, BC V6T 1Z2\\
Canada \\
\printead{e1}\\
\phantom{E-mail: }\printead*{e2}} 
\end{aug}

\thankstext{t1}{Supported in part by NSERC (Canada).}

\thankstext{t2}{Supported in part by the Peter Wall
Institute of Advanced Studies (UBC).}

\received{\smonth{9} \syear{2009}}
\revised{\smonth{2} \syear{2010}}

%
\begin{abstract}
Let $M_n$ be the number of steps of the loop-erasure of a simple random
walk on $\mathbb{Z}^2$ from the origin to the circle of radius $n$. We
relate the moments of $M_n$ to $\operatorname{Es}(n)$, the probability
that a random
walk and an independent loop-erased random walk both started at the
origin do not intersect up to leaving the ball of radius $n$. This
allows us to show that there exists $C$ such that for all $n$ and all
$k=1,2,\ldots, \mathbf{E} [ M_n^k ] \leq C^k k!
\mathbf{E} [ M_n ]^k$ and hence to
establish exponential moment bounds for $M_n$. This implies that there
exists $c>0$ such that for all $n$ and all $\lambda\geq0$,
\[
\mathbf{P} \{ M_n > \lambda\mathbf{E} [ M_n
] \} \leq2e^{-c \lambda}.
\]
Using similar techniques, we then establish a second moment result for
a specific conditioned random walk which enables us to prove that for
any $\alpha< 4/5$, there exist $C$ and $c'>0$ such that for all $n$
and $\lambda> 0$,
\[
\mathbf{P} \{ M_n < \lambda^{-1} \mathbf{E} [ M_n
] \} \leq C e^{-c' \lambda
^{\alpha}}.
\]
\end{abstract}

%

\end{frontmatter}

\section{Introduction}\label{sec1}

The loop-erased random walk (LERW) is a process obtained by
chronologically erasing loops from a random walk on a graph. Since its
introduction by Lawler \cite{Law80}, this process has played a
prominent role in the statistical physics literature. It is closely
related to other models in statistical physics and, in particular, to
the uniform
spanning tree (UST).
Pemantle \cite{Pem91} proved that the unique path
between any two vertices $u$ and $v$ on the UST has the same
distribution as a LERW from $u$ to $v$ and Wilson \cite{Wil96}
devised a powerful algorithm to construct the UST using LERWs.
The existence of a scaling limit of LERW on $\mathbb{Z}^d$ is now known
for all $d$.
For $d \geq4$, Lawler \cite{Law91,Law95} showed that LERW
scales to Brownian motion.
For $d=2$, Lawler, Schramm and Werner \cite{LSW04} proved
that LERW has a conformally invariant scaling limit, Schramm--Loewner evolution; indeed,
LERW was the prototype for the definition of SLE by Schramm \cite{Sch00}.
Most recently, for $d=3$, Kozma \cite{Koz07} proved that
the scaling limit exists and is invariant under rotations and
dilations.

Let $S[0,\sigma_n]$ be simple random walk on $\mathbb{Z}^2$ started
at the
origin and stopped
at~$\sigma_n$, the first time $S$ exits $B_n$, the ball of radius $n$
with center
the origin.
Let $M_n$ be the number of steps of $\mathrm{L}(S[0,\sigma
_n])$, the loop-erasure
of $S[0,\sigma_n]$. In \cite{Ken00}, using domino tilings, Kenyon
proved, for simple random walk on $\mathbb{Z}^2$, that
%
%
\begin{equation} \label{exponent}
\lim_{n \to\infty} \frac{\log\mathbf{E} [ M_n
]}{\log n} = \frac{5}{4}.
\end{equation}
Using quite different methods, Masson \cite{Mas08} extended this
to irreducible bounded
symmetric random walks on any discrete lattice of $\mathbb{R}^2$.
The quantity $5/4$ is called the \textit{growth exponent} for planar
loop-erased random walk.
We remark that while $\operatorname{SLE}_2$ has Hausdorff dimension
$5/4$ almost surely (see \cite{Beff08}), there is no direct proof
of (\ref{exponent}) from this fact; however, unlike the arguments in
\cite{Ken00},
the approach of
\cite{Mas08} does use the connection between the LERW and
$\operatorname{SLE}_2$.

In this paper, we will not be concerned with the exact value of
$\mathbf{E} [ M_n ]$,
but rather with the obtaining of tail bounds on $M_n$.
Our results hold for more general sets than balls.
Let $D$ be a domain in $\mathbb{Z}^2$ with $D \neq\varnothing,\mathbb
{Z}^2$. Write
$S[0,\sigma_D]$ for simple random walk run until its first exit from $D$,
$\mathrm{L}(S[0,\sigma_D])$ for its loop erasure and
$M_D$ for the number of steps in $\mathrm{L}(S[0,\sigma_D])$.
\begin{theorem} \label{main1}
There exists $c_0>0$ such that the following holds.
Let $D$ be a simply connected subset of $\mathbb{Z}^2$ containing $0$ such
that for
all $z \in D$, $\operatorname{dist}(z,D^c) \leq n$. Then:
\begin{enumerate}
\item
%
%
\begin{equation} \label{eq:meanexp}
\mathbf{E} \bigl[ e^{c_0 M_D/\mathbf{E} [ M_n ]}
\bigr] \leq2;
\end{equation}
\item
consequently, for all $\lambda\geq0$,
%
%
\begin{equation} \label{eq:etail}
\mathbf{P} \{ M_D > \lambda\mathbf{E} [ M_n
] \} \leq2 e^{-c_0 \lambda}.
\end{equation}
\end{enumerate}
\end{theorem}
\begin{theorem} \label{main2}
For all $\alpha< 4/5$, there exist
$C_1(\alpha) <\infty$, $c_2(\alpha)>0$
such that for all $\lambda> 0$, all $n$ and all $D \supset B_n$,
%
%
\begin{equation}\label{eq:etail2}
\mathbf{P} \{ M_D < \lambda^{-1} \mathbf{E} [ M_n
] \} \leq C_1(\alpha)
\exp( -c_2(\alpha) \lambda^{\alpha}).
\end{equation}
\end{theorem}

These results are proven in Theorems \ref{mainupper} and \ref{main},
where a slightly more general situation is considered.
\begin{remarks}
1. We expect that these results will hold for irreducible random walks
with bounded, symmetric increments on any discrete lattice of
$\mathbb{R}^2$. Almost all of the proofs in this paper can be extended to
this more general case without any modification. The one exception is
Lemma \ref{reflection}, where we use the fact that simple random walk on
$\mathbb{Z}^2$ is invariant under reflections with respect to
horizontal and
vertical lines.
Theorem \ref{main1} does not depend on Lemma \ref{reflection} and
therefore should be valid in this generality.
It is likely that an alternative proof of Lemma \ref{reflection}
could be found, but we do not pursue this point further here,
restricting our attention to simple random walk on $\mathbb{Z}^2$.

2. The bound (\ref{eq:etail2}) for general $D \supset B_n$ does not
follow immediately
from (\ref{eq:etail2}) for $B_n$.
The reason is that if $Y$ is $\mathrm{L}(S[0, \sigma_D])$ run
until its first
exit from $B_n$, then $Y$ does not, in general, have the same law
as $\mathrm{L}(S[0, \sigma_n])$. Similar considerations apply
to Theorem \ref{main1}.

3. We also have similar bounds for the infinite loop-erased walk;
see Theorems \ref{mainupper} and \ref{main}.

4. One motivation for proving these results for general domains in
$\mathbb{Z}^2$, rather than just balls,
is to study the uniform spanning tree (UST) via Wilson's algorithm.
In particular, we are interested in the volume of balls
in the intrinsic metric on the UST and this requires estimating the
number of steps of an LERW until
it hits the boundary of a fairly general domain in $\mathbb{Z}^2$.
\end{remarks}

For the remainder of this Introduction, we discuss the case where
$D=B_n$.
The proofs of Theorems \ref{main1} and \ref{main2} involve estimates
of the
higher moments of $M_n$. Building on \cite{Mas08},
we relate $\mathbf{E} [ M_n^k ]$ to $\operatorname
{Es}(n)$, the probability that an LERW and
an independent random walk do not intersect up to leaving the ball of
radius $n$.
We show that there exists $C < \infty$ such that
%
%
\begin{eqnarray} \label{eq:intro3}
\mathbf{E} [ M_n^k ] &\leq& C^k k! (n^2 \operatorname
{Es}(n))^k \qquad\mbox{(Theorem \ref
{kmoment});} \\
\label{eq:intro4}
\mathbf{E} [ M_n ] &\geq& C n^2
\operatorname{Es}(n) \qquad\mbox
{(Proposition \ref{M_Dlb}).}
\end{eqnarray}
It is not surprising that the moments of $M_n$ are related to
$\operatorname{Es}
(n)$. To begin with,
\[
\mathbf{E} [ M_n^k ]
= \sum_{z_1, \ldots, z_k \in B_n} \mathbf{P} \{ z_1, \ldots
, z_k \in\mathrm{L} (S[0,\sigma_n]) \}.
\]
Furthemore, for a point $z$ to be on $\mathrm{L}(S[0,\sigma
_n])$, it must be
on the random walk
path $S[0,\sigma_n]$ and not be on the loops that get erased. In order
for this to occur,
the random walk path after $z$ cannot intersect the loop-erasure of the
random walk path
up to $z$. Therefore, for $z$ to be on $\mathrm{L}(S[0,\sigma
_n])$, a random
walk and an
independent LERW must not intersect in a neighborhood of $z$.
Generalizing this
to $k$ points, we get for each $i$, a contribution of $\operatorname
{Es}(r_i)$, where
$r_i$ is chosen small enough to give ``near independence'' of events in the
balls $B_{r_i}(z_i)$.
Propositions \ref{kformula} and\vspace*{2pt} \ref{kpoints} make this
approach precise.
Summing over the $C^k n^{2k}$ $k$-tuples of points in $B_n$ and using facts
about $\operatorname{Es}(\cdot)$ that we establish in Section
\ref{esprobsec} gives
(\ref{eq:intro3}).

Combining (\ref{eq:intro3}) and (\ref{eq:intro4}) yields
\[
\mathbf{E} [ M_n^k ] \leq C^k k! \mathbf{E}
[ M_n ]^k \qquad\mbox{(Theorem \ref
{mainupper}),}
\]
from which Theorem \ref{main1} follows easily.

To establish (\ref{eq:etail2}), we prove a second moment bound for a specific
conditioned random walk and combine this with an iteration argument, as follows.
Let $B_n(x)$ be the ball of radius $n$ centered at $x \in\mathbb
{Z}^2$ and
$R_n$ be the
square $\{(x,y) \in\mathbb{Z}^2\dvtx-n \leq x,y \leq n \}$. Fix a positive
integer $k$ and
consider $\mathrm{L}(S[0,\sigma_{kn}])$. We first establish an
upper bound for
\[
\mathbf{P} \{ M_{kn} < \mathbf{E} [ M_n ]
\}.
\]

Let $k' = k/\sqrt2$ (so that $R_{k'n} \subset B_{kn}$). Let $\gamma
_j$ be
the restriction of $\mathrm{L}(S[0,\sigma_{kn}])$ from $0$ up
to the first exit
of $R_{jn}$, $j=0, \ldots, k'$. For $j=0, \ldots, k'-1$, let $x_j \in
\partial R_{jn}$
be the point where $\gamma_j$ hits $\partial R_{jn}$ and $B_j = B_{n}(x_j)$.
Finally, for $j=1, \ldots, k'$, let $N_j$ be the number of steps of
$\gamma_j$
from $x_{j-1}$ up to the first time it exits $B_{j-1}$;
see Figure \ref{iteratepic}
in Section \ref{lowersec}.
We consider squares instead of balls to take advantage of the symmetry of
simple random walk on $\mathbb{Z}^2$ with respect to vertical and horizontal
lines, as
mentioned above.
Clearly,
%
%
\begin{eqnarray}\label{afr}
\mathbf{P} \{ M_{kn} < \mathbf{E} [ M_n ]
\} &\leq& \mathbf{P} \Biggl( \bigcap_{j=1}^{k'}\{ N_j <
\mathbf{E} [ M_n ] \} \Biggr) \nonumber\\[-8pt]\\[-8pt]
&\leq& \prod_{j=1}^{k'} \max_{\gamma_{j-1}} \mathbf{P} \{
N_j < \mathbf{E} [ M_n ] \mid\gamma_{j-1} \}.\nonumber
\end{eqnarray}
However, by the domain Markov property for LERW (Lemma \ref{condit}),
conditioned
on $\gamma_{j-1}$, the rest of the LERW curve is obtained by running a
random walk
conditioned to leave $B_{kn}$ before hitting $\gamma_{j-1}$ and then
erasing loops.
For this reason, we will be interested in the number of steps of the
loop-erasure
of a random walk started on the boundary of a square and conditioned to
leave some large
ball before hitting a set contained in the square. Formally, we give
the following
definition (throughout this paper, we identify $\mathbb{R}^2$ with
$\mathbb{C}$ and use complex notation such as ``$\arg$'' and
``$\operatorname{Re}$'').
\begin{defin}[(See Figure \ref{definitionpic})] \label{definition}
Suppose that the natural numbers $m$, $n$ and $N$ are such that $\sqrt
{2}m+n \leq N$
and that $K$ is a subset of the square $R_m = [-m,m]^2$. Suppose that
$x = (m,y)$ with $ | y | \leq m$ is any point on
the right-hand side of $R_m$, and let $X$ be a random walk started at $x$,
conditioned to leave $B_N$ before hitting~$K$. Let
%
%
\begin{figure}

\includegraphics{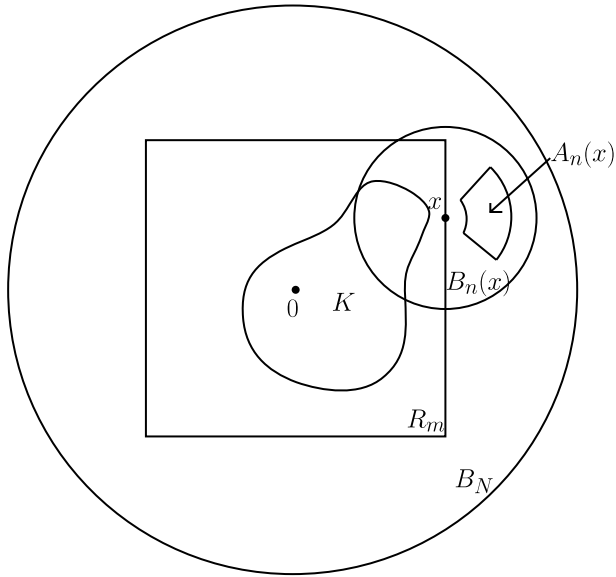}

\caption{Setup for Definition \protect\ref{definition}.}
\label{definitionpic}
\end{figure}
$\alpha$ be $\mathrm{L}(X[0,\sigma_N])$ from $x$ up to its
first exit time of
the ball $B_n(x)$. We then let $M_{m,n,N,x}^{K}$ be the number of steps
of $\alpha$ in
$A_n(x) = \{z\dvtx n/4 \leq| z-x | \leq3n/4, | {\arg}
(z-x) | \leq
\pi/4 \}$.
Note that the condition $\sqrt{2}m + n \leq N$ ensures that $B_n(x)$
is contained in $B_N$.
\end{defin}

We look at the number of steps of the LERW in $A_n(x)$ rather than
in $B_n(x)$ since the expectations of these random variables are
comparable and
it is convenient not to have to worry about points that are close to
$x$, $K$ or $\partial B_n(x)$.
We are therefore interested in estimating
\[
\mathbf{P} \{ M_{m,n,N,x}^{K} < \mathbf{E} [ M_n
] \}.
\]
To do this, we first show that (up to a $\log$ term) $\mathbf{E}
[ M_{m,n,N,x}^{K} ]$ is
comparable to $n^2 \operatorname{Es}(n)$ and, therefore, by (\ref{eq:intro4}),
$\mathbf{E} [ M_{m,n,N,x}^{K} ]$ is comparable to
$\mathbf{E} [ M_n ]$ (Proposition
\ref{ExpN}).
Next, we prove that $\mathbf{E} [ (M_{m,n,N,x}^{K})^2
]$ is comparable
to $\mathbf{E} [ M_{m,n,N,x}^{K} ]^2$ (again up to a
$\log$ term; see
Proposition \ref{ExpN^2}).
By a standard second moment technique, this implies that there exist
$c=c(n,N) > 0$
and $p=p(n,N) < 1$ such that
%
%
\begin{equation}\label{eq:MKEM}
\mathbf{P} \{ M_{m,n,N,x}^{K} < c \mathbf{E} [
M_{m,n,N,x}^{K} ] \} < p.
\end{equation}
Using the fact that $\mathbf{E} [ M_n ]$ is comparable
to $\mathbf{E} [ M_{m,n,N,x}^{K} ]$, we can
then plug this into (\ref{afr}) to conclude that there exists $p=p(k)=
1 -c (\log k)^{-8}$ such that
%
%
\begin{equation} \label{eq:introkn}
\mathbf{P} \{ M_{kn} < \mathbf{E} [ M_n ]
\} < p^k.
\end{equation}
Finally, to prove (\ref{eq:etail2}), one makes an
appropriate choice of $k$ and relates $\mathbf{E} [ M_{kn}
]$ to $\mathbf{E} [ M_n ]$.
Although the logarithmic corrections in Propositions \ref{ExpN} and
\ref{ExpN^2}
mean that $p$ in (\ref{eq:MKEM}) depends on $n$ and $N$, and so $p$ in
(\ref{eq:introkn}) depends on $k$, this correction is small enough so that
(\ref{eq:introkn}) still gives a useful bound.

The paper is organized as follows. In Section \ref{rwdefsec}, we fix
notation and recall the basic
properties of random walks that will be needed. In Section
\ref{lerwsec}, we give a precise definition
of the LERW and state some of its properties. Many of these properties
were established in
\cite{Mas08}. Indeed, this paper uses similar techniques to those in
\cite{Mas08}, most notably, relating
the growth exponent to $\operatorname{Es}(n)$. It turns out that the
latter quantity
is often easier to analyze
directly; see Section \ref{esprobsec}.

Section \ref{greensec} contains some technical lemmas involving
estimates for Green's
functions for random walks in various domains and for the conditioned
random walks $X$ in
Definition \ref{definition}. In Section \ref{expmomentsec}, we prove
Theorem \ref{main1}
using the approach described above.
Finally, in Section \ref{lowersec}, we use the iteration outlined
above to prove
Theorem \ref{main2}.

\section{Definitions and background for random walks}
\label{rwdefsec}

\subsection{Notation for random walks and Markov chains}

Throughout the paper, when we say random walk, we will mean simple
random walk on $\mathbb{Z}^2$. We will denote a random walk starting
at a
point $z \in\mathbb{Z}^2$ by $S^z$. When $z=0$, we will omit the
superscript. If we have two random walks $S^z$ and $S^w$, starting at
two different points $z$ and $w$, then we assume that they are independent
unless otherwise specified. We use similar notation for other Markov
chains on $\mathbb{Z}^2$ (all our Markov chains are assumed to be
time-homogeneous). When there is no possibility of confusion, we
will also use the following standard notation: given an event $A$
that depends on a Markov chain $X$, we let $\mathbf{P}^{z} ( A
)$ denote the
probability of $A$ given that $X_0 = z$.

\subsection{A note about constants}

For the entirety of the paper, we will use the letters $c$ and $C$ to
denote positive constants that
will not depend on any variable, but may change from appearance to appearance.
When we wish to fix a constant, we will number it with a subscript
(e.g., $c_0$).

Given two positive functions $f(n)$ and $g(n)$, we write $f(n) \asymp
g(n)$ if
there exists $C < \infty$ such that for all $n$,
\[
C^{-1} g(n) \leq f(n) \leq C g(n).
\]
We will say that two sequences of events $\{E_n\}$ and $\{F_n\}$ have
the same
probability ``up to constants'' if $\mathbf{P} ( E_n )
\asymp\mathbf{P} ( F_n )$
and are independent ``up to constants'' if $\mathbf{P} ( E_n
\cap F_n )
\asymp\mathbf{P} ( E_n ) \mathbf{P} ( F_n
)$.
We will also use the obvious generalization for two sequences
of random variables to have the same distribution ``up to constants''
and to be independent ``up to constants.''

\subsection{Subsets of $\mathbb{Z}^2$}
\label{subsetsec}

Given two points $x,y \in\mathbb{Z}^2$, we write $x \sim y$ if
\mbox{$| x-y | = 1$.}

A sequence of points $\omega= [\omega_0, \ldots, \omega_k] \subset
\mathbb{Z}^2$ is called a \textit{path} if $\omega_{j-1} \sim\omega
_{j}$ for
$j=1,\ldots,k$. We let $ | \omega| = k$ be the length of
the path,
$\Theta_k$ be the set of paths of length $k$ and $\Theta= \bigcup_k
\Theta_k$ denote the set of all finite paths. Also, if $X$ is a Markov
chain with transition probabilities $p^X(\cdot,\cdot)$ and $\omega\in
\Theta
_k$, then we define
\[
p^X(\omega) = \prod_{i=1}^k p^X(\omega_{i-1},\omega_i).
\]
Thus, if $X = S$ is a simple random walk, $p^S(\omega) = 4^{-k}$. A
set $D \subset\mathbb{Z}^2$ is connected if, for any pair of points
$x, y
\in
D$, there exists a path $\omega\subset D$ connecting $x$ and~$y$, and
$D$ is simply connected if it is connected and all connected
components of $\mathbb{Z}^2\setminus D$ are infinite.

Given $z \in\mathbb{Z}^2$, let
\[
B_n(z) = \{ x \in\mathbb{Z}^2\dvtx | x - z | \leq n \}
\]
be the ball of radius $n$ centered at $z$ in $\mathbb{Z}^2$. We will
write $B_n$
for $B_n(0)$ and sometimes write $B(z;n)$ for $B_n(z)$. Also, let $R_n$ denote
the square $\{(x,y) \in\mathbb{Z}^2\dvtx-n \leq x,y \leq n \}$.

The outer boundary of a set $D \subset\mathbb{Z}^2$ is
\[
\partial D = \{ x \in\mathbb{Z}^2\setminus D \mbox{: there exists
$y \in
D$ such that $x \sim y$} \}
\]
and its inner boundary is
\[
\partial_i D = \{ x \in D \mbox{: there exists $y \in\mathbb{Z}^2
\setminus D$ such that $x \sim y$} \}.
\]
We also write $\overline{D} = D \cup\partial D$.

Given a Markov chain $X$ on $\mathbb{Z}^2$ and a set $D \subset
\mathbb{Z}^2$, let
\[
\sigma^X_D = \min\{j \geq1 \dvtx X_j \in\mathbb{Z}^2\setminus D \}
\]
be the first exit time of the set $D$ and
\[
\xi^X_D = \min\{j \geq1 \dvtx X_j \in D \}
\]
be the first hitting time of the set $D$. We let $\sigma^X_n = \sigma^X_{B_n}$
and use a similar convention for $\xi^X_n$. If $X$ is
a random walk $S^z$ starting at $z \in\mathbb{Z}^2$, then we let
$\sigma^z_D$
and $\xi^z_D$ be the exit and hitting times for $S^z$. If $z = 0$,
then we
will omit the superscripts. We will also omit superscripts when it is
clear what process the stopping times refer to. For instance, we will
write $X[0, \sigma_n]$ instead of $X[0, \sigma^X_n]$.

\subsection{Basic facts about random walks}

For a Markov chain $X$ and $x,y \in D \subset\mathbb{Z}^2$, let
\[
G_D^X(x,y) = \mathbf{E}^{x} \Biggl[ \sum_{j=0}^{\sigma^X_D - 1}
\mathbh{1}\{X_j = y \} \Biggr]
\]
denote Green's function for $X$ in $D$. We will sometimes write
$G^X(x,y; D)$ for $G^X_D(x,y)$.
We will write $G^X_n(x,y)$ for $G^X_{B_n}(x,y)$ and when $X = S$ is a random
walk, we will omit the superscript $S$.

Recall that a function $f$ defined on $\overline{D} \subset\mathbb
{Z}^2$ is
discrete harmonic on $D$
if, for all $z \in D$,
\[
\mathcal{L}f(z) := -f(z) + \frac{1}{4} \sum_{x \sim z} f(x) = 0.
\]
For any two disjoint subsets $K_1$ and $K_2$ of $\mathbb{Z}^2$, it is
easy to
verify that the function
\[
h(z) = \mathbf{P}^{z} \{ \xi_{K_1} < \xi_{K_2} \}
\]
is discrete harmonic on $\mathbb{Z}^2\setminus(K_1 \cup K_2)$.
Furthermore, if we let $X$ be a random walk conditioned to hit $K_1$
before $K_2$, then $X$
is a reversible Markov chain on $\mathbb{Z}^2\setminus(K_1 \cup K_2)$ with
transition probabilities
\[
p^X(x,y) = \frac{1}{4} \frac{h(y)}{h(x)}.
\]
Therefore, if $\omega= [\omega_0, \ldots, \omega_k]$ is a path in
$\mathbb{Z}^2\setminus(K_1 \cup K_2)$, then
%
%
\begin{equation}\label{htransform}
p^X(\omega) = \frac{h(\omega_k)}{h(\omega_0)} 4^{- | \omega
|}.
\end{equation}
Using this fact, the following lemma follows readily.
\begin{lemma} \label{greencondit} Suppose that $X$ is a random walk
conditioned to
hit $K_1$ before $K_2$ and let $D$ be such that
$D \subset\mathbb{Z}^2\setminus(K_1 \cup K_2)$. Then, for any $x,y
\in D$,
\[
G^X_D(x, y) = \frac{h(y)}{h(x)} G_D(x,y).
\]
In particular, $G^X_D(x,x) = G_D(x,x)$.
\end{lemma}

Using a last-exit decomposition, one can also express $h(x)$ in terms
of Green's functions; see \cite{Mas08}, Lemma 3.1.
\begin{lemma} \label{last-exit}
Let $K_1, K_2 \subset\mathbb{Z}^2$ be disjoint and $x \in\mathbb
{Z}^2\setminus
(K_1 \cup K_2)$. Then,
\begin{eqnarray*}
&& \mathbf{P}^{x} \{ \xi_{K_1} < \xi_{K_2} \} \\
&&\qquad= \frac{G (x, x; \mathbb{Z}^2\setminus(K_1 \cup K_2)
)}{G(x, x;
\mathbb{Z}^2\setminus K_1)}
\sum_{y \in\partial_i K_1} \mathbf{P}^{y} \{ \xi_x < \xi
_{K_2} \mid\xi_x < \xi_{K_1} \} \mathbf{P}^{x} \{
S(\xi_{K_1}) = y \}.
\end{eqnarray*}
\end{lemma}

The following proposition was proven in \cite{Mas08} and will be used
frequently in the paper.
\begin{prop} \label{dirichlet}
There exists $c > 0$ such that for all $n$ and all
$K \subset\{z \in\mathbb{Z}^2\dvtx\operatorname{Re}(z) \leq0 \}$,
\[
\mathbf{P} \biggl\{ \arg(S(\sigma_n)) \in\biggl[-\frac{\pi}{4},
\frac{\pi}{4} \biggr] \Bigm|\sigma_n < \xi_{K} \biggr\} \geq c.
\]
\end{prop}

We conclude this section with a list of standard potential theory
results that will be used
throughout the paper, often without referring back to this proposition.
The proofs of these results can all be found in \cite{LL08}, Chapter 6.
\begin{prop} \label{p:potential}
\begin{enumerate}
\item
(Discrete Harnack principle.) Let $U$ be a connected open subset of
$\mathbb{R}^2$ and
$A$ a compact subset of $U$. There then exists a constant $C(U,A)$ such
that for
all $n$ and all positive harmonic functions $f$ on $nU \cap\mathbb{Z}^2$,
\[
f(x) \leq C(U,A) f(y)
\]
for all $x, y \in nA \cap\mathbb{Z}^2$.
\item
There exists $c>0$ such that for all $n$ and all paths
$\alpha$ connecting $B_n$ to $\mathbb{Z}^2\setminus B_{2n}$,
\begin{eqnarray*}
\mathbf{P}^{z} \{ \xi_\alpha< \sigma_{2n} \} &\geq& c
\qquad\mbox{for all $z
\in B_n$,} \\
\mathbf{P}^{z} \{ \xi_\alpha< \xi_{n} \} &\geq& c
\qquad\mbox{for all $z \in
\partial B_{2n}$.}
\end{eqnarray*}
\item
If $m < | z | < n$, then
\[
\mathbf{P}^{z} \{ \xi_m < \sigma_n \} = \frac{\ln n -
{\ln}| z | +
O(m^{-1})}{\ln n - \ln m}.
\]
\item
If $z \in B_n$, then
\[
\mathbf{P}^{z} \{ \xi_0 < \sigma_n \} = \biggl(1 -
\frac{{\ln}| z |}{\ln n}
\biggr)
\biggl[1 + O \biggl(\frac{1}{\ln n} \biggr) \biggr].
\]
\item
If $z \in B_n \setminus\{0\}$, then
\[
G_n(0,z) \asymp\ln\frac{n}{ | z |}.
\]
\item
\[
G_n(0,0) \asymp\ln n.
\]
\end{enumerate}
\end{prop}

\section{Loop-erased random walks}
\label{lerwsec}

\subsection{Definition}

We now describe the loop-erasing procedure and define the loop-erased
random walk. Given a
path $\lambda= [\lambda_0, \ldots, \lambda_m]$ in $\mathbb{Z}^2$, let
$\mathrm{L}(\lambda) = [\widehat{\lambda}_0, \ldots, \widehat
{\lambda}_{n}]$ denote its
chronological loop-erasure. More precisely, we let
\[
s_0 = \sup\{ j \dvtx\lambda(j) = \lambda(0) \}
\]
and, for $i > 0$,
\[
s_i = \sup\{ j\dvtx\lambda(j) = \lambda(s_{i-1} + 1) \}.
\]
Let
\[
n = \inf\{ i\dvtx s_i = m \}.
\]
Then,
\[
\mathrm{L} (\lambda) = [\lambda(s_0), \lambda(s_1), \ldots,
\lambda(s_n)].
\]

One may obtain a different result if one performs the loop-erasing
procedure backward instead of forward. In other words, if we let
$\lambda^R = [\lambda_m, \ldots, \lambda_0]$ be the time reversal of
$\lambda$, then, in general,
\[
\mathrm{L}^R(\lambda) := (\mathrm{L}(\lambda
^R) )^R \neq\mathrm{L}(\lambda).
\]
However, the following lemma shows that if $\lambda$ is distributed
according to a Markov chain, then
$\mathrm{L}^R(\lambda)$ has the same distribution as
$\mathrm{L}(\lambda)$. Recall that
$\Theta$ denotes the set of all finite paths in $\mathbb{Z}^2$.
\begin{lemma}[(Lawler \cite{Law91})] \label{lerwbf} There exists a
bijection $T \dvtx\Theta\to\Theta$ such that
\[
\mathrm{L}^R(\lambda) = \mathrm{L}(T \lambda).
\]
Furthermore, $T \lambda$ and $\lambda$ visit the same edges in
$\mathbb{Z}^2
$ in the same directions with the same
multiplicities so that, for any Markov chain $X$ on $\mathbb{Z}^2$, $p^X(T
\lambda) = p^X(\lambda)$.
\end{lemma}

A fundamental fact about LERWs is the following ``domain Markov property.''
\begin{lemma}[(Domain Markov property \cite{Law91})] \label{condit}
Let $D \subset\Lambda$ and $\omega= [\omega_0, \omega_1,\break \ldots,
\omega_k]$ be a path in $D$.
Let $Y$ be a random walk started at $\omega_k$ conditioned to exit $D$
before hitting $\omega$.
Suppose that $\omega' = [\omega'_0 ,\ldots, \omega'_{k'}]$ is such that
\[
\omega\oplus\omega' := [\omega_0, \ldots, \omega_k, \omega'_0,
\ldots, \omega'_{k'}]
\]
is a path from $\omega_0$ to $\partial D$. Then, if we let $\alpha$
be the first $k$ steps of $\mathrm{L}(S[0,\sigma_D])$,
\[
\mathbf{P} \{ \mathrm{L}(S[0,\sigma_D]) = \omega\oplus
\omega'
\mid\alpha= \omega\}
= \mathbf{P} \{ \mathrm{L}(Y[0,\sigma_D]) = \omega'
\}.
\]
\end{lemma}

Suppose that $l$ is a positive integer and $D$ is a proper subset of
$\mathbb{Z}^2$ with \mbox{$B_l \subset D$}. Let $\Omega_l$ be the set
of paths
$\omega= [0, \omega_1, \ldots, \omega_k] \subset\mathbb{Z}^2$
such that
$\omega_j \in B_l$, $j=1, \ldots, k-1$, and $\omega_k \in\partial
B_l$. Define the
measure $\mu_{l,D}$ on $\Omega_l$ to be the distribution on $\Omega
_l$ obtained by restricting
$\mathrm{L}(S[0,\sigma_D])$ to the part of the path from $0$ to
the first
exit of $B_l$.

Two different sets $D_1$ and $D_2$ will produce different measures.
However, the following
proposition \cite{Mas08} shows that as $\mathbb{Z}^2\setminus D_1$ and
$\mathbb{Z}^2\setminus D_2$
get farther away from~$B_l$, the measures $\mu_{l,D_1}$ and $\mu
_{l,D_2}$ approach each other.
\begin{prop}\label{ax}
There exists $C < \infty$ such that the following holds.
Suppose that $n \geq4$, $D_1$ and $D_2$ are such that
$B_{nl} \subset D_1$ and $B_{nl} \subset D_2$, and $\omega\in\Omega
_l$. Then,
\[
1 - \frac{C}{\log n} \leq\frac{\mu_{l, D_1}(\omega)}{\mu
_{l,D_2}(\omega)}
\leq1 + \frac{C}{\log n}.
\]
\end{prop}

The previous proposition shows that for a fixed $l$, the sequence
$\mu_{l,n}(\omega) := \mu_{l,B_n}(\omega)$ is Cauchy. Therefore,
there exists a
limiting measure $\mu_l$ such that
\[
\lim_{n \to\infty} \mu_{l,n}(\omega) = \mu_l(\omega).
\]
The $\mu_l$ are consistent and therefore there exists a measure $\mu$ on
infinite self-avoiding paths. We call the associated process the infinite
LERW and denote it by $\widehat{S}$. We denote the exit time of a set $D$
for $\widehat{S}$ by $\widehat{\sigma}_D$. An immediate corollary of
the previous
proposition and the definition of $\widehat{S}$ is the following.
\begin{cor} \label{infdist} Suppose that $B_{4l} \subset D$ and
$\omega\in\Omega_l$.
Then,
\[
\mathbf{P} \{ \widehat{S}[0,\widehat{\sigma}_l] = \omega
\} \asymp\mu
_{l,D}(\omega).
\]
\end{cor}

The following result follows immediately from Corollary
\ref{infdist} and \cite{Mas08}, Proposition 4.2.
\begin{cor} \label{condinf}
Suppose that $B_{4l} \subset D_1$ and $B_{4l} \subset D_2$, and let $X$
be a random
walk conditioned to leave $D_1$ before $D_2$.
Let $\alpha$ be $\mathrm{L}(X[0,\sigma_{D_1}])$
from $0$ up to its first exit of $B_l$. Then, for $\omega\in\Omega_l$,
\[
\mathbf{P} \{ \alpha= \omega\} \asymp\mathbf
{P} \{ \widehat{S}[0,\widehat{\sigma}_l] = \omega\}.
\]
\end{cor}

We conclude this section with a ``separation lemma'' for random walks
and LERWs. It states the intuitive fact that, conditioned on the event
that a random walk $S$ and an independent infinite LERW $\widehat{S}$ do
not intersect up to leaving $B_n$, the probability that they are
farther than some fixed distance apart from each other on $\partial
B_n$ is bounded from below by $p > 0$.
\begin{prop}[(Separation lemma \cite{Mas08})] \label{sep}
There exist $c, p > 0$ such that for all $n$, the following holds.
Let $S$ and $\widehat{S}$ be independent and let
\[
d_n = \operatorname{dist}(S(\sigma_n), \widehat{S}[0, \widehat
{\sigma}_n])
\wedge\operatorname{dist}(\widehat{S}(\widehat{\sigma}_n), S[0,
\sigma_n]).
\]
Then,
%
%
\begin{equation} \label{eq:sep}
\mathbf{P} \{ d_n \geq cn \mid S[1, \sigma_n] \cap\widehat
{S}[0, \widehat{\sigma}_n] = \varnothing\} \geq p.
\end{equation}
\end{prop}

\subsection{Escape probabilities for LERW}
\label{esprobsec}

\begin{defin}
For a set\vspace*{2pt} $D$ containing $0$, we let $M_D$ be the number
of steps
of $\mathrm{L}(S[0,\sigma_D])$ and $M_n = M_{B_n}$.
We also\vspace*{1pt} let $\widehat{M}_D$ be the number of steps of
$\widehat
{S}[0,\widehat
{\sigma}_D]$
and $ \widehat{M}_n = \widehat{M}_{B_n}$.
\end{defin}

As described in the \hyperref[sec1]{Introduction}, one of the goals of
this paper is to
relate the
moments of $M_D$ and $\widehat{M}_D$ to escape probabilities, which we
now define.
\begin{defin}
Let $S$ and $S'$ be two independent random walks started at $0$. For $m
\leq n$,
let $ \mathrm{L}(S'[0,\sigma_n]) = \eta= [0, \eta_1, \ldots,
\eta_k]$,
$k_0 = \max\{j \geq1 \dvtx\eta_j \in B_m \}$ and $\eta_{m,n}(S') =
[\eta_{k_0}, \ldots, \eta_{k}]$. We then define
\begin{eqnarray*}
\operatorname{Es}(m,n) &=& \mathbf{P} \{ S[1,\sigma_n] \cap
\eta_{m,n}(S') = \varnothing\},\\
\operatorname{Es}(n) &=& \mathbf{P} \{ S[1,\sigma_n] \cap
\mathrm{L}(S'[0,\sigma_n]) = \varnothing\} ,\\
\operatorname{\widehat{Es}}(n) &=& \mathbf{P} \{ S[1,\sigma
_n] \cap\widehat{S}[0, \widehat{\sigma}_n] = \varnothing\}.
\end{eqnarray*}
We also let $\operatorname{Es}(0) = 1$.
\end{defin}

Thus, $\operatorname{Es}(m,n)$ is the probability that a random walk
from the origin
to $\partial B_n$
and the terminal part of an independent LERW from $m$ to $n$ do not intersect.
$\operatorname{Es}(n)$ is the probability that a random walk from the
origin to
$\partial B_n$
and the loop-erasure of an independent random walk from the origin
to $\partial B_n$ do not intersect. $\operatorname{\widehat{Es}}(n)$
is the
corresponding escape probability
for an \textit{infinite} LERW from the origin to $\partial B_n$.

The following was proven in \cite{Mas08}; see Lemma 5.1, Propositions
5.2, 5.3 and Theorem 5.6.
\begin{theorem} \label{bigthm} There exists $C < \infty$ such that
the following hold:
\begin{enumerate}
\item
\[
C^{-1} \operatorname{Es}(n) \le\operatorname{\widehat{Es}}(n) \le C
\operatorname{Es}(n);
\]
\item
for all $k \geq1$, there exists $N = N(k)$ such that for $n \geq N$,
\[
C^{-1}k^{-3/4} \leq\operatorname{Es}(n, kn) \leq C k^{-3/4};
\]
\item
for all $l \leq m \leq n$,
\[
C^{-1} \operatorname{Es}(n) \le\operatorname{Es}(m) \operatorname
{Es}(m,n) \le C \operatorname{Es}(n)
\]
and
\[
C^{-1} \operatorname{Es}(l,n) \le\operatorname{Es}(l,m)
\operatorname{Es}(m,n) \le C \operatorname{Es}(l,n).
\]
\end{enumerate}
\end{theorem}

We conclude this section with some easy consequences of this theorem.
\begin{lemma} \label{Es(kn)} For all $k \geq1$, there exists $c(k) >
0$ such that for all $n = 1, 2, \ldots,$
\[
\operatorname{Es}(kn) \geq c(k) \operatorname{Es}(n).
\]
\end{lemma}
\begin{pf} By parts 2 and 3 of Theorem \ref
{bigthm}, there exists $N(k)$ such that for $n \geq N(k)$,
\[
\operatorname{Es}(kn) \geq c \operatorname{Es}(n,kn) \operatorname
{Es}(n) \geq c k^{-3/4} \operatorname{Es}(n) = c(k) \operatorname{Es}(n).
\]
Since there are only finitely many $n \leq N(k)$, the result holds.
\end{pf}
\begin{lemma} \label{Esdecreasing} There exists $C < \infty$ such
that for all $l \leq m \leq n$,
%
%
\begin{equation} \label{Esdecreasing1}
\operatorname{Es}(n) \leq C \operatorname{Es}(m)
\end{equation}
and
%
%
\begin{equation} \label{Esdecreasing2}
\operatorname{Es}(l,n) \leq C \operatorname{Es}(l,m).
\end{equation}
\end{lemma}
\begin{pf}
Using Theorem \ref{bigthm}, part 3 and the fact that
$\operatorname{Es}(m,n) \leq1$, one obtains that
\[
\operatorname{Es}(n) \leq C \operatorname{Es}(m) \operatorname
{Es}(m,n) \leq C \operatorname{Es}(m)
\]
and
\[
\operatorname{Es}(l,n) \leq C \operatorname{Es}(l,m) \operatorname
{Es}(m,n) \leq C \operatorname{Es}(l,m).
\]
\upqed
\end{pf}
\begin{lemma} \label{agrt} For all $\varepsilon> 0$, there exist
$C(\varepsilon) < \infty$ and $N(\varepsilon)$
such that for all $N(\varepsilon) \leq m \leq n$,
\[
C(\varepsilon)^{-1} \biggl( \frac{n}{m} \biggr)^{-3/4 - \varepsilon}
\leq\operatorname{Es}(m,n) \leq C(\varepsilon) \biggl( \frac{n}{m}
\biggr)^{-3/4 +
\varepsilon}.
\]
\end{lemma}
\begin{pf}
Fix $\varepsilon> 0$. Let $C_1$ be the largest of the constants in the
statements of Theorem \ref{bigthm} and
Lemma \ref{Esdecreasing} and let $j$ be any integer greater
than $C_1^{2/\varepsilon}$. By Theorem \ref{bigthm}, part 2,
there exists $N$ such that for all $n \geq N$,
\[
C_1^{-1} j^{-3/4} \leq\operatorname{Es}(n, jn) \le C_1 j^{-3/4}.
\]
We will show that the conclusion of the lemma holds with this choice of $N$.

Let $m$ and $n$ be such that $N \leq m \leq n$ and let $k$ be the
unique integer such that
\[
j^k \leq\frac{n}{m} < j^{k+1}.
\]
It follows from Theorem \ref{bigthm}, part 3 and Lemma
\ref{Esdecreasing} that
\begin{eqnarray*}
\operatorname{Es}(m,n) &\leq& C \operatorname{Es}(m, j^{k}m) \\
&\leq& C_1^{k+1} \prod_{i=0}^{k-1} \operatorname{Es}(j^{i}m, j^{i+1}m)
\\
&\leq& C_1^{2k+1} (j^{-3/4})^k \\
&\leq& C_1 j^{\varepsilon k} j^{3/4} \biggl( \frac{n}{m}
\biggr)^{-3/4} \\
&\le& C_1 j^{3/4} \biggl( \frac{n}{m} \biggr)^{\varepsilon} \biggl(
\frac
{n}{m} \biggr)^{-3/4}.
\end{eqnarray*}
This proves the upper bound with $C(\varepsilon) = C_1 j^{3/4}$; the
lower bound is proved in exactly the same way.
\end{pf}
\begin{lemma} \label{nEs(n)}
For all $\varepsilon> 0$, there exists $C(\varepsilon) < \infty$
such that
for all $m \leq n$,
\[
m^{3/4+ \varepsilon} \operatorname{Es}(m) \leq C(\varepsilon) n^{3/4
+ \varepsilon} \operatorname{Es}(n).
\]
\end{lemma}
\begin{pf}
Fix $\varepsilon> 0$. Applying Lemma \ref{agrt}, we get that there exist
$c > 0$ and $N$ such that for all $N \leq m \leq n$,
\[
\operatorname{Es}(m,n) \geq c \biggl( \frac{n}{m}
\biggr)^{-3/4-\varepsilon}.
\]
Therefore, if $N \leq m \leq n$, then, by Theorem \ref{bigthm}, part
3,
\begin{eqnarray*}
n^{3/4+\varepsilon} \operatorname{Es}(n) &\geq& c n^{3/4 +
\varepsilon} \operatorname{Es}(m) \operatorname{Es}(m,n) \geq c
n^{3/4+\varepsilon} \operatorname{Es}(m) \biggl(
\frac{n}{m} \biggr)^{-3/4-\varepsilon}\\
&=& c m^{3/4+\varepsilon} \operatorname{Es}(m).
\end{eqnarray*}
Since there are only finitely many pairs $(m,n)$ such that $m \leq n
\leq N$, there
exists $C$ such that $m^{3/4+\varepsilon} \operatorname{Es}(m) \leq C
n^{3/4+\varepsilon
} \operatorname{Es}(n)$ for all such
pairs $(m,n)$. Finally, if $m \leq N \leq n$, then, since
$m^{3/4+\varepsilon} \operatorname{Es}(m) \leq C N^{3/4+\varepsilon}
\operatorname{Es}(N)$ and\break $
N^{3/4+\varepsilon} \operatorname{Es}(N) \leq C n^{3/4+\varepsilon}
\operatorname{Es}(n)$,
the result also holds in this case.
\end{pf}

In Sections \ref{expmomentsec} and \ref{lowersec}, we will have to
handle various sums involving $\operatorname{Es}(n)$
and we will use the following result many times.
\begin{cor} \label{Es_sums}
Let $\gamma> 0$, $\beta>0$ and $1+\alpha- 3\gamma/4 >0$.
There then exists $C < \infty$ (depending on $\alpha$, $\beta$,
$\gamma$) such that for all $n$,
\[
\sum_{j=1}^n j^\alpha\biggl(\ln\frac{n}{j} \biggr)^\beta
\operatorname{Es}
(j)^\gamma
\le C n^{\alpha+1} \operatorname{Es}(n)^\gamma.
\]
\end{cor}
\begin{pf}
Choose $\varepsilon>0$ such that $1+\alpha- 3\gamma/4 - (\beta+
\gamma)
\varepsilon>0$.
Then, using Lemma \ref{nEs(n)},
\begin{eqnarray*}
\sum_{j=1}^n j^\alpha\biggl(\ln\frac{n}{j} \biggr)^\beta
\operatorname{Es}
(j)^\gamma
&=& \sum_{j=1}^n j^{\alpha- 3\gamma/4- \gamma\varepsilon}
\biggl(\ln\frac
{n}{j} \biggr)^\beta(j^{3/4 + \varepsilon} \operatorname
{Es}(j) )^\gamma\\
&\le& C (n^{3/4 + \varepsilon} \operatorname{Es}(n)
)^\gamma\sum_{j=1}^n
j^{\alpha- 3\gamma/4- \gamma\varepsilon} ({n}/{j})^{\varepsilon
\beta} \\
&\le& C n^{3\gamma/4 + \varepsilon\gamma+\varepsilon\beta}
\operatorname{Es}(n)^\gamma
\sum_{j=1}^n j^{\alpha- 3\gamma/4- \gamma\varepsilon-\beta
\varepsilon} \\
&\le& C n^{1+\alpha} \operatorname{Es}(n)^\gamma.
\end{eqnarray*}
\upqed
\end{pf}

\section{Green's function estimates}
\label{greensec}

\begin{lemma} \label{l:greensum}
There exists $C < \infty$ such that the following holds.
Let $D \subset\mathbb{Z}^2$ and, for $z \in D$,
write $\operatorname{dist}(z, D^c)$ for the distance between $z$ and
$D^c$. Let
\[
D_n =\{ z \in D\dvtx\operatorname{dist}(z, D^c) \leq n\}.
\]
Suppose that for all $z \in D_n$, there exists a path in $D^c$
connecting $B(z,n+1)$ to $B(z,2n)^c$.
Then, for any $w \in D$,
%
%
\begin{equation}\label{e:GDsum}
\sum_{z \in D_n} G_D(w,z) \le C n^2.
\end{equation}
In particular, if $D$ is simply connected and $\operatorname{dist}(z,
D^c) \leq n$
for all $z \in D$, then, for all $w \in D$,
\[
\sum_{z \in D} G_D(w,z) \leq C n^2.
\]
\end{lemma}
\begin{pf} Fix $n \geq1$ and define stopping times $(T_j)$, $(U_j)$ as follows:
%
%
\begin{eqnarray}\label{e:tj+1}
T_1 &=& \min\{ i \ge0\dvtx S_i \in D_n \}; \nonumber\\
U_j &=& \min\{ i \ge T_j\dvtx|S_i - S_{T_j}| \ge2 n \}; \\
T_{j+1} &=& \min\{U_j \leq i < \sigma_D \dvtx S_i \in D_n \}.\nonumber
\end{eqnarray}
Here, as usual, we take $T_{j+1}=\infty$ if the set in (\ref{e:tj+1})
is empty.
On the event that $T_j < \infty$, $\mathbf{E}^{S_{T_j}} [ U_j
-T_j ] \le C
n^2$ and thus
%
%
\begin{eqnarray}\label{e:GDp-bnd}
\sum_{z \in D_n} G_D(w,z) &=& \mathbf{E}^{w} \Biggl[ \sum
_{j=1}^{\sigma_D-1} \mathbh{1} \{X_j \in D_n \} \Biggr]\nonumber\\
&\le&\mathbf{E}^{w} \Biggl[ \sum_{j=1}^\infty(U_j- T_j ) \Biggr] \\
&\le& C n^2 \sum_{j=1}^\infty\mathbf{P}^{w} \{ T_j< \infty
\}.\nonumber
\end{eqnarray}
By Proposition \ref{p:potential}, part 2
and our assumption that for all $z \in D_n$, there is a path in $D^c$ connecting
$B(z, n+1)$ and $B(z, 2n)^c$, there exists $p>0$ such that
for any $z \in D_n$,
\[
\mathbf{P}^{z} \bigl\{ \sigma_D < \sigma_{B(z,2n)} \bigr\} > p.
\]
Consequently, $ \mathbf{P}^{w} \{ T_{j+1} < \infty\mid T_j <
\infty\} < 1-p$
and so $\mathbf{P}^{w} \{ T_j< \infty\} < (1-p)^{j-1}$.
Therefore, summing the
series in (\ref{e:GDp-bnd})
yields (\ref{e:GDsum}).
\end{pf}
\begin{lemma} \label{l:probleaveD}
There exist $C < \infty$ and $c>0$ such that the following holds.
Suppose that $D \subset\mathbb{Z}^2$, $w \in D$ is such\vadjust
{\goodbreak} that
$\operatorname{dist}(w,D^c) = n$ and there exists a path in $D^c$ connecting
$B(w,n+1)$ to $B(w,2n)^c$. Then:
\begin{enumerate}
\item
for all $z \in B_{n/2}(w)$,
%
%
\begin{equation} \label{eq:probleaveD}
\mathbf{P}^{z} \{ \xi_{w} < \sigma_D \} \leq C \mathbf
{P}^{z} \bigl\{ \xi_{w} < \sigma_D \wedge\sigma_{B_{2n}(w)}
\bigr\};
\end{equation}
\item
for all $z \in B_n(w)$ and $l \leq| z-w |$,
%
%
\begin{equation} \label{eq:probleaveD5}
\mathbf{P}^{z} \bigl\{ \sigma_D < \xi_{B_l(w)} \bigr\} \geq c
\mathbf{P}^{z} \bigl\{ \sigma_D \wedge\sigma_{B_{2n}(w)} < \xi
_{B_l(w)} \bigr\}.
\end{equation}
\end{enumerate}
\end{lemma}
\begin{pf}
We can take $w = 0$ so that $\sigma_{B_{2n}(w)} = \sigma_{2n}$ and
$\xi_{B_l(w)} = \xi_l$. We begin with (\ref{eq:probleaveD}).
Let $z_0 \in\partial B_{n/2}$ be such that
\[
\mathbf{P}^{z_0} \{ \xi_{0} < \sigma_D \}
= \max_{z \in\partial B_{n/2}} \mathbf{P}^{z} \{ \xi_{0} <
\sigma_D \}.
\]
Then,
%
%
\begin{eqnarray}\label{eq:probleaveD3}
\mathbf{P}^{z_0} \{ \xi_{0} < \sigma_D \}
&\leq&\mathbf{P}^{z_0} \{ \xi_{0} < \sigma_D \wedge\sigma
_{2n} \} \nonumber\\
&&{}+ \mathbf{P}^{z_0} \{ \sigma_{2n} < \sigma_D
\} \max_{y \in\partial B_{2n}}
\mathbf{P}^{y} \{ \xi_{0} < \sigma_D \} \\
&\leq&\mathbf{P}^{z_0} \{ \xi_{0} < \sigma_D \wedge\sigma
_{2n} \} + \mathbf{P}^{z_0} \{ \sigma_{2n} < \sigma_D
\} \mathbf{P}^{z_0} \{ \xi_0 < \sigma_D \}.\nonumber
\end{eqnarray}
By our assumption, there exists a path in $D^c$ connecting $\partial
B_n$ to $\partial B_{2n}$ and therefore, by Proposition \ref
{p:potential}, part 2, there exists $c > 0$ such that
\[
\mathbf{P}^{z_0} \{ \sigma_{2n} < \sigma_D \} \leq1 - c.
\]
Thus, inserting this in (\ref{eq:probleaveD3}) yields
\[
\mathbf{P}^{z_0} \{ \xi_{0} < \sigma_D \} \leq C
\mathbf{P}^{z_0} \{ \xi_{0} < \sigma_D \wedge\sigma_{2n}
\}.
\]
Hence, if $z$ is any point in $\partial B_{n/2}$, we have
%
%
\begin{eqnarray} \label{eq:probleaveD4}
\mathbf{P}^{z} \{ \xi_0 < \sigma_D \} &\leq& \mathbf
{P}^{z_0} \{ \xi_{0} < \sigma_D \}
\leq C \mathbf{P}^{z_0} \{ \xi_{0} < \sigma_D \wedge\sigma
_{2n} \} \nonumber\\[-8pt]\\[-8pt]
&\leq& C \mathbf{P}^{z} \{ \xi_{0} < \sigma_D
\wedge\sigma_{2n} \},\nonumber
\end{eqnarray}
where the last inequality follows from the discrete Harnack inequality.

Now, suppose that $z$ is any point in $B_{n/2}$. Then, using (\ref
{eq:probleaveD4}),
we have
\begin{eqnarray*}
\mathbf{P}^{z} \{ \xi_0 < \sigma_D \} &=& \mathbf
{P}^{z} \{ \xi_0 < \sigma_{n/2} \}
+ \sum_{y \in\partial B_{n/2}} \mathbf{P}^{y} \{ \xi_0 <
\sigma_D \} \mathbf{P}^{z} \{ \sigma_{n/2} < \xi_0;
S(\sigma_{n/2}) = y \} \\
&\leq& \mathbf{P}^{z} \{ \xi_0 < \sigma_{n/2} \}\\
&&{} + C \sum_{y \in\partial B_{n/2}} \mathbf{P}^{y} \{ \xi_0 <
\sigma_D \wedge\sigma_{2n} \} \mathbf{P}^{z} \{
\sigma_{n/2} < \xi_0; S(\sigma_{n/2}) = y \} \\
&\leq& C \biggl( \mathbf{P}^{z} \{ \xi_0 < \sigma_{n/2}
\} \\
&&\hspace*{12pt}{} + \sum_{y \in
\partial B_{n/2}} \mathbf{P}^{y} \{ \xi_0 < \sigma_D \wedge
\sigma_{2n} \} \mathbf{P}^{z} \{ \sigma_{n/2} < \xi_0;
S(\sigma_{n/2}) = y \} \biggr) \\
&=& C \mathbf{P}^{z} \{ \xi_0 < \sigma_D \wedge\sigma_{2n}
\}.
\end{eqnarray*}
This proves (\ref{eq:probleaveD}).

The proof of (\ref{eq:probleaveD5}) is simpler. By Proposition \ref
{p:potential}, part 2,
\begin{eqnarray*}
\mathbf{P}^{z} \{ \sigma_D < \xi_l \} &\geq& \mathbf
{P}^{z} \{ \sigma_D < \xi_l \wedge\sigma_{2n} \} +
\mathbf{P}^{z} \{ \sigma_{2n} < \xi_l \wedge\sigma_D
\} \min
_{y \in\partial B_{2n}} \mathbf{P}^{y} \{ \sigma_D < \xi_n
\} \\
&\geq& \mathbf{P}^{z} \{ \sigma_D < \xi_l \wedge\sigma_{2n}
\} + c \mathbf{P}^{z} \{ \sigma_{2n} < \xi_l \wedge
\sigma_D \} \\
&\geq& c ( \mathbf{P}^{z} \{ \sigma_D < \xi_l \wedge
\sigma_{2n} \} + \mathbf{P}^{z} \{ \sigma_{2n} < \xi_l
\wedge\sigma_D \} ) \\
&=& c \mathbf{P}^{z} \{ \sigma_D \wedge\sigma_{2n} < \xi_l
\}.
\end{eqnarray*}
\upqed
\end{pf}
\begin{lemma} \label{l:greenD}
There exists $C < \infty$ such that the following holds.
Suppose that $D \subset\mathbb{Z}^2$, $w \in D$ is such that
$\operatorname{dist}(w,D^c) = n$ and there exists a path in $D^c$ connecting
$B(w,n+1)$ to $B(w,2n)^c$.
Then, for any $z \in B_{n/2}(w)$,
%
%
\begin{equation} \label{eq:greenD}
G_D(w,z) \leq C G_{B_{2n}(w) \cap D}(w,z).
\end{equation}
\end{lemma}
\begin{pf}
This follows immediately from Lemma \ref{l:probleaveD} and the facts that
\[
G_D(w,z) = \mathbf{P}^{z} \{ \xi_w < \sigma_D \}
\mathbf{P}^{w} \{ \sigma_D < \xi_w \}^{-1}
\]
and
\[
G_{B_{2n}(w) \cap D}(w,z) = \mathbf{P}^{z} \bigl\{ \xi_w < \sigma_D
\wedge\sigma_{B_{2n}(w)} \bigr\} \mathbf{P}^{w} \bigl\{ \sigma
_D \wedge\sigma_{B_{2n}(w)} < \xi_w \bigr\}^{-1}.
\]
\upqed
\end{pf}

Given $D \subset\mathbb{Z}^2$, let $D_{+} = \{z \in D \dvtx
\operatorname
{Re}(z) > 0 \}$ and
$D_{-} = \{z \in D \dvtx\operatorname{Re}(z) < 0 \}$. If $z = (z_1, z_2)
\in\mathbb{Z}^2$, then we
let $\overline{z} = (-z_1, z_2)$ be the reflection of $z$ with respect
to the $y$-axis $\ell$
and $\overline{D} = \{ \overline{z} \dvtx z \in D\}$ be the reflection of
the set $D$.
\begin{lemma}[(See Figure \ref{reflect})] \label{reflection}
Suppose that $K \subset D \subset\mathbb{Z}^2$ are such that $ D_{+}
\subset
\overline{D}_{-}$
%
%
\begin{figure}

\includegraphics{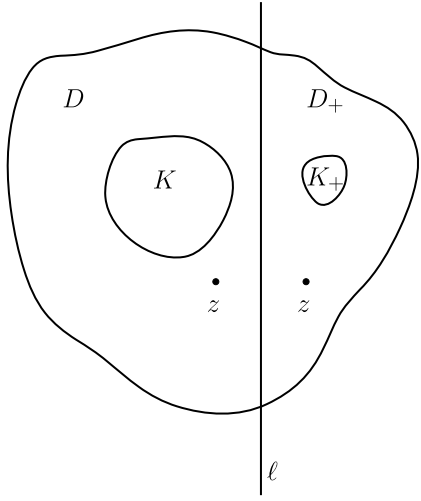}

\caption{The setup for Lemma \protect\ref{reflection}.}
\label{reflect}
\end{figure}
and $K_{+} \subset\overline{K}_{-}$. Then, for all $z \in D_{-}$,
\[
\mathbf{P}^{z} \{ \sigma_D < \xi_K \} \leq\mathbf
{P}^{\overline{z}} \{ \sigma_D < \xi_K \}.
\]
\end{lemma}
\begin{pf}
The proof uses a simple reflection argument. For a random walk started
at $z \in D_{-}$ to escape $D$ before hitting $K$, either it escapes
$D$ before hitting $K$ while staying to the left of $\ell$ or it hits
$\ell$ before hitting $K$ and then escapes $D$ before hitting $K$. In
the first case, the reflected random walk path will be a random path
starting at $\overline{z}$, escaping $D$ before hitting $K$. In the
second case, the reflection of the path up to the first time it hits
$\ell$ will avoid $K$ and hit $\ell$ at the same point. By the Markov
property, the distribution of the paths after this point will be the
same.

More precisely, using the fact that the reflection of a simple random
walk across $\ell$ is again a simple random walk, it follows that for
$z \in D_-$,
\[
\mathbf{P}^{z} \{ \sigma_D < \xi_K \} = \mathbf
{P}^{\overline{z}} \{ \sigma_{\overline{D}} < \xi_{\overline
{K}} \}.
\]
However, since $ D_{+} \subset\overline{D}_{-}$ and $K_{+} \subset
\overline{K}_{-}$,
we have
\begin{eqnarray*}
\mathbf{P}^{\overline{z}} \{ \sigma_{\overline{D}} < \xi
_{\overline{K}} \}
&=& \sum_{x \in\partial D_+} \mathbf{P}^{x} \{ \sigma
_{\overline{D}} < \xi_{\overline{K}} \}
\mathbf{P}^{\overline{z}} \{ \sigma_{D_+} < \xi_{\overline
{K}_-} \wedge\xi_\ell; S(\sigma_{D+}) = x \} \\
&&{} + \sum_{y \in\ell}\mathbf{P}^{y} \{ \sigma
_{\overline{D}} < \xi_{\overline{K}} \}
\mathbf{P}^{\overline{z}} \{ \xi_\ell< \xi_{\overline{K}_-}
\wedge\sigma_{D_+} ; S(\xi_{\ell}) = y \} \\
&\leq&\mathbf{P}^{\overline{z}} \{ \sigma_{D_+} < \xi_{K_+}
\wedge\xi_\ell\} \\
&&{} + \sum_{y \in\ell}\mathbf{P}^{y} \{ \sigma_{D} <
\xi_{K} \}
\mathbf{P}^{\overline{z}} \{ \xi_\ell< \xi_{K_+} \wedge
\sigma_{D_+} ; S(\xi_{\ell}) = y \} \\
&=& \mathbf{P}^{\overline{z}} \{ \sigma_{D} < \xi_{K} \}.
\end{eqnarray*}
\upqed
\end{pf}
\begin{cor} \label{reflectcor} There exists $C < \infty$ such that
the following holds.
Suppose that $m$, $n$, $N$, $K$ and $x$ are as in Definition \ref{definition}.
Then, for all $z \in A_n(x)$,
\[
\max_{w \in\partial B_{n/8}(x)} \mathbf{P}^{w} \{ \sigma_N <
\xi_K \} \leq C
\mathbf{P}^{z} \{ \sigma_N < \xi_K \}.
\]
\end{cor}
\begin{pf}
We apply Lemma \ref{reflection} with $\ell= \{(m,k)\dvtx k \in\mathbb{Z}
\}$ replacing the $y$-axis to conclude that
\[
\max_{w \in\partial B_{n/8}(x)} \mathbf{P}^{w} \{ \sigma_N <
\xi_K \}
= \mathop{\max_{w \in\partial B_{n/8}(x)}}_{\operatorname{Re}(w)
\geq m} \mathbf{P}^{w} \{ \sigma_N < \xi_K \}.
\]

If $x = (m,x_2)$, let
\[
D_n(x) = \{ (w_1,w_2) \in\mathbb{Z}^2\dvtx n/16 \leq w_1 - m \leq n/8,
| w_2 - x_2 | \leq n/8 \}.
\]
Then, by again applying Lemma \ref{reflection}, this time with $\ell=
\{(m+n/16,k) \dvtx k \in\mathbb{Z}\}$,
\[
\mathop{\max_{w \in\partial B_{n/8}(x)}}_{\operatorname{Re}(w) \geq
m} \mathbf{P}^{w} \{ \sigma_N < \xi_K \}
\leq\max_{w \in D_n(x)} \mathbf{P}^{w} \{ \sigma_N < \xi_K
\}.
\]
However, by the discrete Harnack inequality, there exists $C < \infty$
such that for all $z \in A_n(x)$ and all $w \in D_n(x)$,
\[
\mathbf{P}^{w} \{ \sigma_N < \xi_K \} \leq C \mathbf
{P}^{z} \{ \sigma_N < \xi_K \}.
\]
\upqed
\end{pf}
\begin{lemma} \label{markovgreens} There exists $C < \infty$ such
that the following holds.
Suppose that $m$, $n$, $N$, $K$, $x$ and $X$ are as in Definition \ref
{definition}.
Then, for any $z \in A_n(x)$,
\[
C^{-1} \leq G_N^X(x,z) \leq C \ln\frac{N}{n}.
\]
\end{lemma}
\begin{pf}
By Lemma \ref{greencondit},
%
%
\begin{eqnarray} \label{markova}
G_N^X(x,z) &=& G_{B_N \setminus K}(x,z) \frac{\mathbf{P}^{z} \{
\sigma_N < \xi_K \}}{\mathbf{P}^{x} \{ \sigma_N < \xi
_K \}} \\
&=& G_{B_N \setminus K}(z,z) \frac{\mathbf{P}^{x} \{ \xi_z <
\sigma_N \wedge\xi_K \}
\mathbf{P}^{z} \{ \sigma_N < \xi_K \}}{\mathbf{P}^{x}
\{ \sigma_N < \xi_K \}}.\nonumber
\end{eqnarray}

To begin with,
%
%
\begin{equation} \label{markovb}
\ln n \asymp G_{B_{n/8}(z)}(z,z) \leq G_{B_N \setminus K}(z,z) \leq
G_{B_{2N}(z)}(z,z) \asymp\ln N.
\end{equation}

Next,
\begin{eqnarray*}
&& \mathbf{P}^{x} \{ \xi_z < \sigma_N \wedge\xi_K \} \\
&&\qquad= \sum_{y \in\partial_i B_{n/8}(z)} \mathbf{P}^{y} \{ \xi_z
< \sigma_N \wedge\xi_K \}
\mathbf{P}^{x} \bigl\{ S\bigl(\xi_{B_{n/8}(z)}\bigr) = y; \xi
_{B_{n/8}(z)} <
\sigma_N \wedge\xi_K \bigr\}.
\end{eqnarray*}
Furthermore, for any $y \in\partial_i B_{n/8}(z)$,
\[
\mathbf{P}^{y} \{ \xi_z < \sigma_N \wedge\xi_K \}
\leq\mathbf{P}^{y} \bigl\{ \xi_z < \sigma_{B_{2N}(z)} \bigr\}
\leq C \frac{\ln(N/n)}{\ln N}
\]
and
\[
\mathbf{P}^{y} \{ \xi_z < \sigma_N \wedge\xi_K \}
\geq\mathbf{P}^{y} \bigl\{ \xi_z < \sigma_{B_{n/4}(z)} \bigr\}
\geq\frac{c}{\ln n}.
\]

Thus,
%
%
\begin{equation} \label{markovc} \frac{c}{\ln n}
\leq\frac{\mathbf{P}^{x} \{ \xi_z < \sigma_N \wedge\xi_K
\}}{\mathbf{P}^{x} \{ \xi_{B_{n/8}(z)} < \sigma_N
\wedge\xi_K \}} \leq C \frac{\ln(N/n)}{\ln N}.
\end{equation}

Next, on the one hand,
\[
\mathbf{P}^{x} \{ \sigma_N < \xi_K \} \geq\sum_{y \in
\partial_i B_{n/8}(z)}
\mathbf{P}^{y} \{ \sigma_N < \xi_K \} \mathbf{P}^{x}
\bigl\{ S\bigl(\xi_{B_{n/8}(z)}\bigr) = y; \xi_{B_{n/8}(z)} < \sigma_N
\wedge\xi_K \bigr\}.
\]
By the discrete Harnack inequality, there exists $C$ such that for any
$y \in\partial_i B_{n/8}(z)$,
\[
\mathbf{P}^{z} \{ \sigma_N < \xi_K \} \leq C \mathbf
{P}^{y} \{ \sigma_N < \xi_K \}.
\]
Therefore,
\[
\mathbf{P}^{x} \{ \sigma_N < \xi_K \} \geq c \mathbf
{P}^{z} \{ \sigma_N < \xi_K \} \mathbf{P}^{x} \bigl\{
\xi_{B_{n/8}(z)} < \sigma_N \wedge\xi_K \bigr\}.
\]
On the other hand,
\[
\mathbf{P}^{x} \{ \sigma_N < \xi_K \} = \sum_{w \in
\partial B_{n/8}(x)}
\mathbf{P}^{w} \{ \sigma_N < \xi_K \} \mathbf{P}^{x}
\bigl\{ S\bigl(\sigma_{B_{n/8}(x)}\bigr) = w; \sigma_{B_{n/8}(x)} < \xi_K
\bigr\}.
\]
By Corollary \ref{reflectcor}, for any $w \in\partial B_{n/8}(x)$,
\[
\mathbf{P}^{w} \{ \sigma_N < \xi_K \} \leq C \mathbf
{P}^{z} \{ \sigma_N < \xi_K \}.
\]
Therefore,
\[
\mathbf{P}^{x} \{ \sigma_N < \xi_K \} \leq C \mathbf
{P}^{z} \{ \sigma_N < \xi_K \} \mathbf{P}^{x} \bigl\{
\sigma_{B_{n/8}(x)} < \xi_K \bigr\}.
\]
Finally, by Proposition \ref{dirichlet},
\begin{eqnarray*}
\mathbf{P}^{x} \bigl\{ \sigma_{B_{n/8}(x)} < \xi_K \bigr\}
&\leq& C \mathbf{P}^{x} \biggl\{ \sigma_{B_{n/8}(x)} < \xi_K; \bigl|
{\arg}\bigl(S\bigl(\sigma_{B_{n/8}(x)}\bigr) - x\bigr) \bigr| \leq\frac
{\pi}{4}
\biggr\} \\
&\leq& C \mathbf{P}^{x} \bigl\{ \xi_{B_{n/8}(z)} < \sigma_N \wedge
\xi_K \bigr\}.
\end{eqnarray*}
Thus,
%
%
\begin{equation} \label{markovd}
\mathbf{P}^{x} \{ \sigma_N < \xi_K \} \asymp\mathbf
{P}^{z} \{ \sigma_N < \xi_K \} \mathbf{P}^{x} \bigl\{
\xi_{B_{n/8}(z)} < \sigma_N \wedge\xi_K \bigr\}.
\end{equation}

The result then follows by combining (\ref{markova}), (\ref
{markovb}), (\ref{markovc}) and (\ref{markovd}).
\end{pf}

\section{Exponential moments for $M_D$ and $\widehat{M}_D$}
\label{expmomentsec}

To reduce the size of our expressions, we use the following notation.
For this section only, we will use the symbol$\mbox{ }\disjoint\mbox{
}$to denote
the disjoint
intersection relation. Thus, if $K_1$ and $K_2$ are two subsets of
$\mathbb{Z}^2$,
we will write $K_1 \disjoint K_2$ to mean $K_1 \cap K_2 = \varnothing$.
\begin{defin} \label{kdefin}
Suppose that $z_0, z_1, \ldots, z_k$ are any distinct points in a
domain $D \subset\mathbb{Z}^2$
and that $X$ is a Markov chain on $\mathbb{Z}^2$ with $\mathbf
{P}^{z_0} \{ \sigma^X_D < \infty\} = 1$.
We then let $E^X_{z_0, \ldots, z_k}$ be the event that $z_1, z_2,
\ldots, z_k$ are all visited by the path $\mathrm{L}(X^{z_0}[0,
\sigma_D])$
in order.
\end{defin}
\begin{prop} \label{kformula} Suppose that $z_0, z_1, \ldots, z_k$
are distinct points in a
domain $D \subset\mathbb{Z}^2$ and $X$ is a Markov chain on $\mathbb
{Z}^2$ with
$\mathbf{P}^{z_0} \{ \sigma^X_D < \infty\} = 1$.
Define $z_{k+1}$ to be
$\partial D$ and for $i=0,\ldots,k$,
let $X^i$ be independent versions of $X$ started at $z_i$ and $Y^i$ be $X^i$
conditioned on the event $\{\xi^{X^i}_{z_{i+1}} \leq\sigma^{X^i}_D \}$.
Let $\tau^i = \max\{l \leq\sigma^{Y^i}_D \dvtx Y^i_l = z_{i+1} \}$. Then,
\[
\mathbf{P} ( E^X_{z_0, \ldots, z_k} )
= \Biggl[ \prod_{i=1}^k G_D^X(z_{i-1},z_i) \Biggr]
\mathbf{P} \Biggl( \bigcap_{i=1}^k \Biggl\{\mathrm
{L}(Y^{i-1}[0,\tau^{i-1}]) \disjoint\bigcup_{j=i}^{k} Y^j[1,\tau^j]
\Biggr\}\Biggr).
\]
\end{prop}
\begin{pf}
We will write the exit times $\sigma^{X^j}_D$ as $\sigma^j_D$ and the
hitting times $\xi^{X^j}_{z_i}$ as $\xi^j_i$, $i,j=0,\ldots,k$.
For $i,j=0,\ldots,k$, we also let
\[
T^j_i = \cases{
\max\{l \leq\sigma^j_D \dvtx X^j_l = z_i \}, &\quad if $\xi^j_{i} <
\sigma^j_D$,\cr
\sigma^j_D, &\quad if $\sigma^j_D \leq\xi^j_{i}$.}
\]
For $i=0, \ldots, k-1$, let
\[
F_i = \{T^i_{i+1} < \cdots< T^i_k < \sigma^i_D\},
\]
and for $i=0, \ldots, k-2$, let
\[
G_i = \bigcap_{j=i+2}^{k} \{ \mathrm{L}(X^i[T^i_{j-1},T^i_j])
\disjoint
X^i(T^i_j,\sigma^i_D] \}.
\]
Then, by the definition of the loop-erasing procedure,
%
%
\begin{eqnarray} \label{eq:kformula}
\mathbf{P} ( E^X_{z_0, \ldots, z_k} ) = \mathbf
{P} \{ F_0; \mathrm{L}(X^0[0,T^0_1]) \disjoint
X^0(T^0_1, \sigma^0_D]; G_0 \}.
\end{eqnarray}

Conditioned on $\{ T^0_1 < \sigma^0_D \}$, $X^0[0, T^0_1]$ and
$X^0[T^0_1,\sigma^0_D]$ are
independent. $X^0[0,T^0_1]$ has the same distribution as $Y^0[0,\tau
^0]$ and
$X^0[T^0_1,\sigma^0_D]$ has the same distribution as $X^1$ conditioned to
leave $D$ before returning to $z_1$.

The event $\{ T^0_1 < \sigma^0_D \}$ is the same as $\{ \xi^0_1 <
\sigma^0_D \}$. Therefore,
\begin{eqnarray*}
\mathbf{P} ( E^X_{z_0, \ldots, z_k} )
&=& \mathbf{P} \{ \xi^0_1 < \sigma^0_D \} \mathbf
{P} \{ F_1; \mathrm{L}
(Y^0[0, \tau^0]) \disjoint X^1[1, \sigma^1_D];\\
&&\hspace*{64.4pt} \mathrm{L}(X^1[0,T^1_2])
\disjoint X^1(T^1_2, \sigma^1_D]; G_1 \mid\sigma^1_D < \xi^1_1
\} \\
&=& \frac{\mathbf{P} \{ \xi^0_1 < \sigma^0_D \}
}{\mathbf{P} \{ \sigma^1_D < \xi^1_1 \}}
\mathbf{P} \{ F_1; \mathrm{L}(Y^0[0, \tau^0])
\disjoint X^1[1, \sigma^1_D];\\
&&\hspace*{65.5pt} \mathrm{L} (X^1[0,T^1_2])
\disjoint X^1(T^1_2, \sigma^1_D]; G_1 \} \\
&=& G^X_D(z_0,z_1) \mathbf{P} \bigl\{ F_1; \mathrm{L}(Y^0[0, \tau
^0]) \disjoint\bigl(X^1[1, T^1_2] \cup X^1(T^1_2,\sigma^1_D]
\bigr) ;\\
&&\hspace*{121pt}{} \mathrm{L}(X^1[0,T^1_2]) \disjoint X^1(T^1_2, \sigma^1_D];
G_1 \bigr\}.
\end{eqnarray*}
By repeating the previous argument $k-1$ times with $X^1, \ldots
,X^{k-1}$, we obtain the desired result.
\end{pf}

Now, suppose that $D' \subset D$ and let $\beta$ be $\mathrm{L}
(X^{z_0}[0,\sigma_D])$ from
$z_0$ up to the first exit time of $D'$. It is possible to generalize
the previous formula to the probability
that $\beta$ hits $z_1,\ldots,z_k$ in order. However, we will only
require this for the case
where $k=1$ and therefore, to avoid introducing any new notation, we
will only state the result in this case.
\begin{lemma} \label{markovformulalb} Suppose that $D' \subset D$, $z$
and $w$ are distinct points in $D'$ and $X$ is a Markov chain started
at $w$.
Suppose, further, that $\mathbf{P}^{w} \{ \sigma_D < \infty
\} = 1$.
Let $Y$ be $X$ conditioned to hit $z$ before leaving $D$ and let $\tau
$ be the last time
that $Y$ visits $z$ before leaving $D$. Then, if $\beta$ is
$\mathrm{L}
(X[0,\sigma_D])$
from $w$ up to the first exit time of $D'$,
\[
\mathbf{P} \{ z \in\beta\} = G_D^X(w,z) \mathbf
{P} \{ \mathrm{L}(Y[0, \tau]) \disjoint X^z[1, \sigma
_D]; \mathrm{L}(Y[0, \tau]) \subset D' \}.
\]
\end{lemma}
\begin{pf}
As in the proof of Proposition \ref{kformula}, let
\[
T_z = \cases{
\max\{l \leq\sigma^X_D \dvtx X_l = z \}, &\quad if $\xi^X_{z} <
\sigma^X_D$,\cr
\sigma^X_D, &\quad if $\sigma^X_D \leq\xi^X_{z}$.}
\]
Then,
\[
\mathbf{P} \{ z \in\beta\} = \mathbf{P}
\{ T_z < \sigma^X_D; \mathrm{L}(X[0, T_z]) \disjoint X[T_z +
1, \sigma_D]; \mathrm{L}(X[0, T_z]) \subset D' \}.
\]
The proof is then identical to that of Proposition \ref{kformula}.
\end{pf}
\begin{defin} \label{r_idefin}
Suppose that $z_0, z_1, \ldots, z_k$ are any points (not necessarily distinct)
in a domain $D \subsetneq\mathbb{Z}^2$ and let $\mathbf{z} = (z_0,
\ldots, z_k)$.
We then define $z_{k+1}$ to be $\partial D$, let $d(z_i) =
\operatorname{dist}(z_i,
D^c)$ and let
\[
r^{\mathbf{z}}_i = d(z_i) \wedge| z_i - z_{i-1} | \wedge
| z_{i} - z_{i+1} |, \qquad i = 1,2,\ldots,k.
\]
In addition, if $\pi$ is an element of the symmetric group $\mathfrak{S}_k$
on $\{1,\ldots,k\}$, then we let $\pi(0) = 0$ and $\pi(\mathbf{z})
= (z_0, z_{\pi(1)}, \ldots, z_{\pi(k)})$.
\end{defin}
\begin{prop} \label{kpoints}
There exists $C < \infty$ such that the following holds. Suppose that either:
\begin{enumerate}
\item$z_0, z_1, \ldots, z_k$ are any points in a domain $D \subsetneq
\mathbb{Z}^2$ and $X$
is a random walk $S$ started at $z_0$; or
\item$m$, $n$, $N$, $K$, $x$ and $X$ are as in Definition \ref{definition},
$z_0 = x$, $D = B_N$ and $z_1, \ldots, z_k$ are in $A_n(x)$.
\end{enumerate}
Then, letting $\mathbf{z} = (z_0, \ldots, z_k)$ and $r^{\mathbf
{z}}_i$ be as in Definition \ref{r_idefin},
%
%
\begin{equation} \label{eq:kpoints} \qquad\mathbf{P} \{ z_1,\ldots
,z_k \in\mathrm{L} (X[0,\sigma_D]) \}
\leq C^k \sum_{\pi\in\mathfrak{S}_k} \prod_{i=1}^k G^X_D\bigl(z_{\pi
(i-1)},z_{\pi(i)}\bigr) \operatorname{Es} \bigl(r^{\pi(\mathbf{z})}_{\pi
(i)} \bigr).
\end{equation}
\end{prop}
\begin{pf}
The proofs of the two cases are almost identical and we will prove them
both at the same time.

\textit{First, suppose that $z_0, \ldots, z_k$ are distinct.}
Recall the definition of $E^X_{z_0, \ldots, z_k}$ from Definition \ref
{kdefin}. Then,
\[
\mathbf{P} \{ z_1,\ldots,z_k \in\mathrm
{L}(X[0,\sigma_D]) \} = \sum_{\pi\in
\mathfrak{S}_k} E^X_{z_0,z_{\pi(1)}, \ldots, z_{\pi(k)}}.
\]
Therefore, if we let $Y^0, \ldots, Y^k$ be as in Proposition
\ref{kformula}, then it suffices to show that
\[
\mathbf{P} \Biggl( \bigcap_{i=1}^k \Biggl\{\mathrm{L}(Y^{i-1}[0,\tau
^{i-1}]) \disjoint\bigcup_{j=i}^{k} Y^j[1,\tau
^j] \Biggr\} \Biggr) \leq C^k \prod_{i=1}^k
\operatorname{Es} (r^{\mathbf{z}}_i ).
\]
For $i=1,\ldots,k$, let $B_i = B(z_i; r^{\mathbf{z}}_i/4)$. Then,
\begin{eqnarray*}
&&
\mathbf{P} \Biggl( \bigcap_{i=1}^k \Biggl\{\mathrm{L}(Y^{i-1}[0,\tau
^{i-1}]) \disjoint\bigcup_{j=i}^{k} Y^j[1,\tau^j]
\Biggr\} \Biggr)\\
&&\qquad\leq\mathbf{P} \Biggl( \bigcap_{i=1}^k \{\mathrm
{L}(Y^{i-1}[0,\tau^{i-1}]) \disjoint Y^i[1,\tau^i] \} \Biggr).
\end{eqnarray*}

Let $T\dvtx\Theta\to\Theta$ be the bijection given in Lemma \ref{lerwbf}.
For all $\lambda\in\Theta$, $p^X(T(\lambda)) = p^X(\lambda)$ and
$T \lambda$
visits the same points as $\lambda$. Thus,
\begin{eqnarray*}
&&\mathbf{P} \Biggl( \bigcap_{i=1}^k \{\mathrm{L}(Y^{i-1}[0,\tau
^{i-1}]) \disjoint Y^i[1,\tau^i] \} \Biggr)
\\
&&\qquad= \mathbf{P} \Biggl( \bigcap_{i=1}^k \{\mathrm{L}(T
\circ Y^{i-1}[0,\tau^{i-1}]) \disjoint(T \circ Y^i[1,\tau^i])
\} \Biggr) \\
&&\qquad= \mathbf{P} \Biggl( \bigcap_{i=1}^k \{\mathrm{L}(Y^{i-1}[0,\tau
^{i-1}]^R) \disjoint Y^i[1,\tau^i] \} \Biggr).
\end{eqnarray*}

For $i=1,\ldots,k$, let $\beta^i$ be the restriction of
$\mathrm{L}
(Y^{i-1}[0,\tau^{i-1}]^R)$ from $z_i$ to the first exit of $B_i$. Then,
\[
\mathbf{P} \Biggl( \bigcap_{i=1}^k \{\mathrm{L}(Y^{i-1}[0,\tau
^{i-1}]^R) \disjoint Y^i[1,\tau^i] \} \Biggr)
\leq\mathbf{P} \Biggl( \bigcap_{i=1}^k \{\beta^i \disjoint
Y^i[1,\sigma_{B_i}] \} \Biggr).
\]
Furthermore, by the domain Markov property (Lemma \ref{condit}),
conditioned on $\beta^i = [\beta^i_0, \ldots, \beta^i_m]$,
$Y^{i-1}[0,\tau^{i-1}]$ is,
in case 1, a random walk started at $z_{i-1}$ and conditioned to hit
$\beta^i_m$ before
$\partial D \cup\{\beta^i_0,\ldots,\beta^i_{m-1}\}$; in case 2, it
is a random walk started at $z_{i-1}$
and conditioned to hit $\beta^i_m$ before $K \cup\partial D \cup\{
\beta^i_0,\ldots,\beta^i_{m-1}\}$.
In either case, by the Harnack principle, $Y^{i-1}[0,\sigma
_{B_{i-1}}]$ and $\beta^i$ are independent ``up to constants'' and thus
\[
\mathbf{P} \Biggl( \bigcap_{i=1}^k \{\beta^i \disjoint
Y^i[1,\sigma_{B_i}] \} \Biggr)
\leq C^k \prod_{i=1}^k \mathbf{P} \{ \beta^i \disjoint
Y^i[1,\sigma_{B_i}] \}.
\]
By another application of the Harnack principle, $Y^{i}[0,\sigma_{B_{i}}]$
has the same distribution, up to constants, as a random walk started at
$z_{i}$ and stopped at its
first exit of~$B_{i}$. Furthermore, by Corollary \ref{infdist}, $\beta
^i$ has the
same distribution, up to constants, as an infinite LERW started at
$z_i$ and stopped at the
first exit of $B_i$. Therefore, for $i=1,\ldots,k$,
\[
\mathbf{P} \{ \beta^i \disjoint Y^i[1,\sigma_{B_i}]
\} \leq C \operatorname{\widehat{Es}
}(r^{\mathbf{z}}_i/4).
\]
Finally, by Theorem \ref{bigthm}, part 1 and
Lemma \ref{Es(kn)}, $\operatorname{\widehat{Es}}(r^{\mathbf
{z}}_i/4) \leq C \operatorname{Es}
(r^{\mathbf{z}}_i)$.

\textit{Now, suppose that $z_0, \ldots, z_k$ are any points in $D$.}
Let
\[
p(\mathbf{z}) = \prod_{i=1}^k G^X_D(z_{i-1},z_{i}) \operatorname
{Es}
(r^{\mathbf{z}}_{i} ).
\]
We will establish (\ref{eq:kpoints}) by induction on $k$.
We have already proven that (\ref{eq:kpoints}) holds for $k=1$.
Now, suppose that (\ref{eq:kpoints}) holds for $k-1$ and suppose that
$z_0, \ldots, z_k$ are not distinct.
Since (\ref{eq:kpoints}) involves a sum over all possible permutations
of the
entries of $\mathbf{z}$, we may assume without loss of generality
that $z_{j} = z_{j+1}$ for some $j$. Let $\mathbf{z^{(j)}}$ be
$\mathbf{z}$ with the $j$th entry
deleted and indexed by $\{0,\ldots,k\} \setminus\{j\}$ (so that $z_i
= z^{(j)}_i$ for all $i \neq j$).
Then, since $r^{\mathbf{z^{(j)}}}_i = r^{\mathbf{z}}_i$ for all $i
\neq j$, $i \neq j+1$,
\[
p(\mathbf{z}) = p\bigl(\mathbf{z^{(j)}}\bigr) \cdot G_D^X(z_j,z_j)
\operatorname{Es} (r^{\mathbf{z}}_j ) \operatorname
{Es} (r^{\mathbf
{z}}_{j+1} ) \operatorname{Es} \bigl(r^{\mathbf
{z^{(j)}}}_{j+1} \bigr)^{-1}.
\]
Since $z_j = z_{j+1}$, we have $r^{\mathbf{z}}_{j} = r^{\mathbf
{z}}_{j+1} = 0$ and, therefore,
$\operatorname{Es}(r^{\mathbf{z}}_j) = \operatorname{Es}(r^{\mathbf
{z}}_{j+1}) = 1$. Also,
$G_D^X(z_j,z_j) \geq1$.
Therefore, $ p(\mathbf{z}) \geq p(\mathbf{z^{(j)}})$.

Now, let $\mathfrak{S}_{A}$ be the symmetric group on the set $A = \{
1, \ldots, k \} \setminus\{j \}$.
There then exists an obvious bijection between $\mathfrak{S}_A$ and
\[
\mathfrak{B} = \{ \pi\in\mathfrak{S}_k\dvtx\pi^{-1}(j+1) = \pi
^{-1}(j) + 1 \}.
\]
Therefore, by our induction hypothesis,
\begin{eqnarray*}
&&
\mathbf{P} \{ z_1, \ldots, z_k \in\mathrm{L}(X[0,\sigma_D]) \} \\
&&\qquad\leq C^{k-1} \sum
_{\pi\in\mathfrak{S}_A} p\bigl(\pi\bigl(\mathbf{z^{(j)}}\bigr)\bigr)
\leq C^k \sum_{\pi\in\mathfrak{S}_A} p\bigl(\pi\bigl(\mathbf
{z^{(j)}}\bigr)\bigr) \\
&&\qquad\leq C^k \sum_{\pi\in\mathfrak{B}} p(\pi(\mathbf{z}))
\leq C^k \sum_{\pi\in\mathfrak{S}_k} p(\pi(\mathbf{z})).
\end{eqnarray*}
\upqed
\end{pf}

Recall that if $D$ is a proper subset of $\mathbb{Z}^2$, then $M_D$
denotes the number of steps of $\mathrm{L}(S[0,\sigma_D])$.
Given $D'
\subset D$, we let $M_{D',D}$ denote the number of steps of
$\mathrm{L}
(S[0,\sigma_D])$ while it is in $D'$ or, equivalently, the number of
points in $D'$ that are on the path $\mathrm{L}(S[0,\sigma_D])$.
\begin{theorem} \label{kmoment} There exists $C < \infty$ such that
the following hold:

1. if we suppose that $D \subset\mathbb{Z}^2$ contains $0$
and $D'
\subset D$ is such that for all $z \in D'$,
there exists a path in $D^c$ connecting $B(z,n+1)$ and $B(z,2n)^c$,
then, for all $k=1,2,\ldots,$
\[
\mathbf{E} [ M_{D',D}^k ] \leq C^k k! (n^2
\operatorname{Es}(n))^k;
\]

2. in particular, if $D$ is simply connected, contains $0$
and, for all $z \in D$, $\operatorname{dist}(z, D^c) \leq n$, then
\[
\mathbf{E} [ M_D^k ] \leq C^k k! (n^2 \operatorname{Es}(n))^k.
\]
\end{theorem}
\begin{pf}
Let $\mathfrak{S}_k$ denote the symmetric group on $k$ elements and
recall the
definition of $r^{\mathbf{z}}_i$ given in Definition \ref{r_idefin}
(here, $z_0 = 0$).
Then, by Proposition \ref{kpoints},
\begin{eqnarray*}
\mathbf{E} [ M_{D',D}^k ] &=& \mathbf{E} \biggl[
\biggl(\sum_{z \in D'} \mathbh{1}\{z \in\mathrm{L}
(S[0,\sigma_D])\} \biggr)^k \biggr] \\
&=& \sum_{z_1 \in D'} \cdots\sum_{z_k \in D'} \mathbf{P} \{
z_1,\ldots,z_k \in\mathrm{L}(S[0,\sigma_D]) \} \\
&\leq& C^k \sum_{\pi\in\mathfrak{S}_k} \sum_{z_1 \in D'} \cdots
\sum_{z_k \in D'}
\prod_{i=1}^k G_D\bigl(z_{\pi(i-1)}, z_{\pi(i)}\bigr) \operatorname
{Es}\bigl(r^{\pi(\mathbf
{z})}_{\pi(i)}\bigr) \\
&=& C^k k! \sum_{z_1 \in D'} \cdots\sum_{z_k \in D'} \prod_{i=1}^k
G_D(z_{i-1}, z_{i}) \operatorname{Es}(r^{\mathbf{z}}_{i}).
\end{eqnarray*}
Therefore, it suffices to show that
%
%
\begin{equation} \label{e:ksum}
\sum_{z_1 \in D'} \cdots\sum_{z_k \in D'} \prod_{i=1}^k
G_D(z_{i-1}, z_{i}) \operatorname{Es}(r^{\mathbf{z}}_{i}) \leq C^k
(n^2 \operatorname{Es}(n))^k.
\end{equation}

Let $f_i = G_D(z_{i-1}, z_{i}) \operatorname{Es}(r^{\mathbf{z}}_{i})$
and $F_j =
\prod_{i=1}^j f_i$. Then, if $d(z) = \operatorname{dist}(z, D^c)$,
we have
%
%
\begin{equation}\label{e:kpts}\qquad
\prod_{i=1}^k G_D(z_{i-1}, z_{i}) \operatorname{Es}(r^{\mathbf{z}}_{i})
=F_{k-1} G_D(z_{k-1},z_{k})
\bigl( \operatorname{Es}\bigl( |z_{k}-z_{k-1}| \wedge d(z_{k} )\bigr)
\bigr).
\end{equation}
Since only the terms $f_k$ and $f_{k-1}$ involve $z_k$, and
$\operatorname{Es}(a \wedge b) \le\operatorname{Es}(a) +
\operatorname{Es}(b)$, we then have
\begin{eqnarray*}
&& \sum_{z_1 \in D'} \cdots\sum_{ z_k \in D'} \prod_{i=1}^k
G_D(z_{i-1}, z_{i}) \operatorname{Es}(r^{\mathbf{z}}_{i}) \\
&&\qquad\le\sum_{z_1 \in D'} \cdots\sum_{z_{k-1} \in D'} F_{k-2}
G_D(z_{k-2},z_{k-1}) \\
&&\qquad\quad\hspace*{66pt}{} \times\sum_{z_k \in D'}
G_D(z_{k-1},z_k)\bigl( \operatorname{Es}\bigl( |z_{k-1}-z_{k-2}| \wedge
d(z_{k-1} ) \bigr)\\
&&\qquad\quad\hspace*{100.17pt}\hspace*{110.68pt}{} +
\operatorname{Es}( |z_{k-1}-z_{k}|)
\bigr) \\
&&\qquad\quad\hspace*{100.17pt}{} \times\bigl( \operatorname{Es}(
|z_{k}-z_{k-1}|) + \operatorname
{Es}(d(z_{k})) \bigr).
\end{eqnarray*}
Multiplying out the final terms in the expression above, we need to
bound the following sums:
%
%
\begin{eqnarray} \label{e:s1}\qquad
S_1 &=& \operatorname{Es}\bigl( |z_{k-1}-z_{k-2}| \wedge d(z_{k-1}
)\bigr)
\sum_{z_k} G_D(z_{k-1},z_k) \operatorname{Es}( |z_{k-1}-z_{k}|), \\
\label{e:s2}
S_2 &=& \operatorname{Es}\bigl( |z_{k-1}-z_{k-2}| \wedge d(z_{k-1} )
\bigr)
\sum_{z_k} G_D(z_{k-1},z_k) \operatorname{Es}( d(z_k) ), \\
\label{e:s3}
S_3 &=& \sum_{z_k} G_D(z_{k-1},z_k) \operatorname{Es}(
|z_{k-1}-z_{k}|)^2, \\
\label{e:s4}
S_4 &=& \sum_{z_k} G_D(z_{k-1},z_k) \operatorname{Es}(
|z_{k-1}-z_{k}|) \operatorname{Es}( d(z_k) ).
\end{eqnarray}
Since $2ab \le a^2 + b^2, $ we can bound $S_4$ by
%
%
\begin{equation} \label{e:s5}
S_4 \le S_3 + \sum_{z_k} G_D(z_{k-1},z_k) \operatorname{Es}( d(z_k)
)^2 = S_3 + S_5.
\end{equation}

We first consider $S_3$. Let $D_1 = D \cap B_{n/2}(z_{k-1})$ and $D_2 =
D' \setminus D_1$. Then,
\[
S_3 \leq\sum_{z_k \in D_1} G_D(z_{k-1},z_k) \operatorname{Es}(
|z_{k-1}-z_{k}|)^2
+ \sum_{z_k \in D_2} G_D(z_{k-1},z_k) \operatorname{Es}( |z_{k-1}-z_{k}|)^2.
\]
However, by our assumptions on $D'$ and $D$, and Lemma \ref{l:greenD},
for all $z_{k} \in D_1$, we have
\[
G_D(z_{k-1},z_k) \le C G_{B_{2n}(z_{k-1})}(z_{k-1},z_k) \le C \ln
\biggl( \frac{2n}{|z_{k-1}-z_k|} \biggr).
\]
So,
\begin{eqnarray*}
&&\sum_{z_k \in D_1} G_D(z_{k-1},z_k) \operatorname{Es}(
|z_{k-1}-z_{k}|)^2\\
&&\qquad\le C \sum_{z_k \in D_1} \ln\biggl( \frac{2n}{|z_{k-1}-w|}
\biggr)
\operatorname{Es}( |z_{k-1}-z_k |)^2 \\
&&\qquad\le C \sum_{z_k \in B_{2n}(z_{k-1})} \ln\biggl( \frac
{2n}{|z_{k-1}-z_k|} \biggr) \operatorname{Es}( |z_{k-1}- z_k |)^2 \\
&&\qquad\le C \sum_{j=1}^{2n} j \ln\biggl( \frac{2n}{j} \biggr)
\operatorname{Es}( j
)^2 \\
&&\qquad\le C n^2 \operatorname{Es}(n)^2,
\end{eqnarray*}
where the\vspace*{1pt} last inequality is justified by Corollary \ref{Es_sums}.
Furthermore, for $z_k \in D_2$, $\operatorname{Es}(|z_{k-1} - z_k|)^2
\leq C \operatorname{Es}
(n)^2$. Therefore, by Lemma \ref{l:greensum},
\[
\sum_{z_k \in D_2} G_D(z_{k-1},z_k) \operatorname{Es}( |z_{k-1}-z_{k}|)^2
\leq C \operatorname{Es}(n)^2 \sum_{z_k \in D'} G_D(z_{k-1},z_k) \leq
C n^2 \operatorname{Es}(n)^2.
\]
Therefore, $S_3 \leq C n^2 \operatorname{Es}(n)^2$.
Similarly, we obtain
%
%
\begin{equation}
S_1 \le C \operatorname{Es}\bigl( |z_{k-1}-z_{k-2}| \wedge d(z_{k-1}
)\bigr) n^2
\operatorname{Es}(n).
\end{equation}

Let $D_j = \{z \in D \dvtx d(z) \leq j\}$ be as in Lemma \ref{l:greensum}.
By first applying Lemma \ref{l:greensum} and then Lemma \ref{nEs(n)},
we then have
\begin{eqnarray*}
S_5 &\leq& \sum_{j=0}^{\lceil\log_2 n \rceil} \sum_{z_k \in D_{2^j}
\setminus D_{2^{j-1}}} G_D(z_{k-1}, z_k) \operatorname{Es}( d(z_k) )^2
\\
&\leq& C \sum_{j=0}^{\lceil\log_2 n \rceil} \operatorname
{Es}(2^j)^2 \sum_{z_k
\in D_{2^j} \setminus D_{2^{j-1}}} G_D(z_{k-1}, z_k) \\
&\le& C \sum_{j=0}^{\lceil\log_2 n \rceil} \operatorname{Es}(2^j)^2
\sum_{z_k
\in D_{2^j}} G_D(z_{k-1}, z_k) \\
&\le& C \sum_{j=0}^{\lceil\log_2 n \rceil} 2^{2j} \operatorname
{Es}( 2^{j})^2 \\
&\le& C \sum_{j=0}^{\lceil\log_2 n \rceil} ( (2^{j})^{3/4
+\varepsilon} \operatorname{Es}( 2^{j}) )^2 ( 2^{j})^{1/2 -
2\varepsilon} \\
&\le& C ( n^{3/4 +\varepsilon} \operatorname{Es}(n) )^2
\sum_{j=1}^{\lceil\log
_2 n \rceil} ( 2^j)^{1/2 - 2\varepsilon} \\
&\le& C n^2 \operatorname{Es}(n)^2 .
\end{eqnarray*}
A similar calculation gives
%
%
\begin{equation}
S_2 \le C \operatorname{Es}\bigl( |z_{k-1}-z_{k-2}| \wedge d(z_{k-1}
)\bigr) n^2
\operatorname{Es}(n).
\end{equation}
Combining these bounds gives
\begin{eqnarray*}
&&\sum_{z_1 \in D'} \cdots\sum_{ z_k \in D'} \prod_{i=1}^k
G_D(z_{i-1}, z_{i}) \operatorname{Es}(r^{\mathbf{z}}_{i}) \\
&&\qquad\leq C n^2 \operatorname{Es}(n) \sum_{z_1 \in D'} \cdots\sum_{
z_{k-1} \in D'} F_{k-2}
G_D(z_{k-2},z_{k-1})\\
&&\qquad\quad\hspace*{113.13pt}{}\times\bigl( \operatorname{Es}\bigl(
|z_{k-1}-z_{k-2}| \wedge d(z_{k-1} )\bigr) +
\operatorname{Es}(n) \bigr) \\
&&\qquad\le C n^2 \operatorname{Es}(n) \sum_{z_1 \in D'} \cdots\sum
_{z_{k-1} \in D'}
F_{k-2} G_D(z_{k-2},z_{k-1})\\
&&\qquad\quad\hspace*{114.1pt}{}\times
\bigl( \operatorname{Es}\bigl( |z_{k-1}-z_{k-2}| \wedge d(z_{k-1} )\bigr
) \bigr).
\end{eqnarray*}
Since this is of the same form as (\ref{e:kpts}), except with only
$k-1$ terms,
iterating this argument gives (\ref{e:ksum}).
\end{pf}
\begin{prop} \label{M_Dlb}
There exists $c > 0$ such that for all $n$ and all simply connected $D
\supset B_n$,
\[
\mathbf{E} [ M_D ] \geq c n^2 \operatorname{Es}(n).
\]
\end{prop}
\begin{pf}
By Lemma \ref{nEs(n)}, $n^2 \operatorname{Es}(n)$ is increasing (up
to a constant).
Therefore, we may assume
that $n$ is the largest integer such that $B_n \subset D$.
Let $A_n = \{z\dvtx n/4 \leq| z | \leq3n/4, | {\arg z}
| \leq\pi
/4 \}$ be as in Definition \ref{definition}.
Then, since there are on the order of $n^2$ points in $A_n$, it
suffices to show that for all $z \in A_n$,
%
%
\begin{equation}
\label{eq:M_Dlb} \mathbf{P} \{ z \in\mathrm{L}(S[0,\sigma_D]) \} \geq c
\operatorname{Es}(n).
\end{equation}

By Proposition \ref{kformula},
%
%
\begin{equation}
\label{eq:M_Dlb2} \mathbf{P} \{ z \in\mathrm{L}(S[0,\sigma_D]) \} =
G_D(0,z) \mathbf{P} \{
\mathrm{L}(Y[0,\tau]) \cap S^z[1,\sigma_D] = \varnothing
\},
\end{equation}
where $Y$ is a random walk started at $0$, conditioned to hit $z$
before leaving $D$
and $\tau= \max\{ k < \sigma_D\dvtx Y_k = z \}$. By Lemma \ref{lerwbf},
$\mathrm{L}(Y[0,\tau])$
has the same distribution as $\mathrm{L}(Y[0,\tau]^R)$.
Furthermore, if we
let $Z$ be a random walk started at $z$,
conditioned to hit $0$ before leaving $D$, then $Y[0,\tau]^R$ has the
same distribution as $Z[0,\xi_0]$. Therefore,
\[
\mathbf{P} \{ \mathrm{L}(Y[0,\tau]) \cap S^z[1,\sigma
_D] = \varnothing\} = \mathbf{P} \{ \mathrm{L}(Z[0,\xi_0]) \cap
S^z[1,\sigma_D] = \varnothing\}.
\]
Furthermore,
\[
G_D(0,z) \geq G_n(0,z) \geq c.
\]
Therefore, in order to show (\ref{eq:M_Dlb}), it is sufficient to
prove that
%
%
\begin{equation} \label{eq:M_Dlb3}
\mathbf{P} \{ \mathrm{L}(Z[0,\xi_0]) \cap
S^z[1,\sigma_D] = \varnothing\} \geq c
\operatorname{Es}(n).
\end{equation}

Let $B = B(z;n/8)$ and let $\beta$ be the restriction of
$\mathrm{L}(Z[0, \xi_0])$
from $z$ up to the first time it leaves the ball $B$. Then,
\begin{eqnarray*}
&& \mathbf{P} \{ \mathrm{L}(Z[0,\xi_0]) \cap
S^z[1,\sigma_D] = \varnothing\} \\
&&\qquad= \mathbf{P} \{ \mathrm{L}(Z[0, \xi_0]) \cap S^z[1,
\sigma_D] = \varnothing\mid\beta\cap S^z[1, \sigma_B] = \varnothing
\} \\
&&\qquad\quad{}\times\mathbf{P} \{ \beta\cap S^z[1, \sigma_B] =
\varnothing\}.
\end{eqnarray*}
By Corollary \ref{infdist}, $\beta$ has the same distribution ``up to
constants'' as an
infinite LERW started at $z$ and stopped at the first exit of $B$.
Therefore, by Theorem \ref{bigthm}, part 1 and Lemma \ref
{Esdecreasing},
\[
\mathbf{P} \{ \beta\cap S^z[1, \sigma_B] = \varnothing
\} \geq c \operatorname{\widehat{Es}
}(n/8) \geq c \operatorname{Es}(n/8) \geq c \operatorname{Es}(n).
\]

By the domain Markov property (Lemma \ref{condit}), if we condition on
$\beta$, the rest
of $\mathrm{L}(Z[0, \xi_0])$ is obtained by running a random
walk conditioned to hit $0$
before $\beta\cup\partial D$ and then loop-erasing. Therefore, by the
separation lemma (Proposition \ref{sep})
and Proposition \ref{dirichlet}, there is a probability greater than
$c > 0$ that this conditioned random
walk reaches $\partial B_{n/16}$ without hitting $S^z[1,\sigma_{n}]$
or leaving $B_{7n/8}$.

Therefore, it remains to show that for all $v \in\partial B_{n/16}$,
%
%
\begin{equation} \label{eq:M_Dlb4}
\mathbf{P}^{v} \{ \xi_0 < \sigma_{B_{n/8}} \mid\xi_0 <
\sigma_D \} \geq c
\end{equation}
and for all $w \in\partial B_n$,
%
%
\begin{equation} \label{eq:M_Dlb5}
\mathbf{P}^{w} \{ \sigma_D < \xi_{B_{7n/8}} \} \geq c.
\end{equation}
By Lemma \ref{l:probleaveD},
\[
\mathbf{P}^{v} \{ \xi_0 < \sigma_D \} \leq C \mathbf
{P}^{v} \{ \xi_0 < \sigma_{2n} \}
\]
and
\[
\mathbf{P}^{w} \{ \sigma_D < \xi_{7n/8} \} \geq c
\mathbf{P}^{w} \{ \sigma_{2n} < \xi_{7n/8} \}.
\]
By Proposition \ref{p:potential}, these imply (\ref{eq:M_Dlb4}) and
(\ref{eq:M_Dlb5}).
\end{pf}

Recall the definitions of $M_D$ and $M_{D',D}$ given before Theorem
\ref{kmoment} and recall that $\widehat{M}_n$ denotes the number of steps
of $\widehat{S}[0,\widehat{\sigma}_n]$.
\begin{theorem} \label{mainupper}
There exist $C_0, C_1 < \infty$ and $c_0, c_1 > 0$ such that the
following holds.
Suppose that $D \subset\mathbb{Z}^2$ contains $0$ and $D' \subset D$
is such
that for all $z \in D'$,
there exists a path in $D^c$ connecting $B(z,n+1)$ and $B(z,2n)^c$.
Then:
\begin{enumerate}
\item
for all $k=1,2,\ldots,$
%
%
\begin{equation} \label{eq:emean}
\mathbf{E} [ M_{D',D}^k ] \leq(C_0)^k k! (\mathbf
{E} [ M_n ])^k;
\end{equation}
\item there exists $c_0>0$ such that
%
%
\begin{equation}\label{eq:emean2}
\mathbf{E} \bigl[ \exp\{c_0M_{D',D}/\mathbf{E} [
M_n ] \} \bigr]\le2;
\end{equation}
\item for all $\lambda\ge0$,
%
%
\begin{equation} \label{eq:emean3}
\mathbf{P} \{ M_{D',D} > \lambda\mathbf{E} [ M_n
] \} \leq2 e^{-c_0 \lambda};
\end{equation}
\item for all $n$ and all $\lambda\geq0$,
%
%
\begin{equation}
\mathbf{P} \{ \widehat{M}_n > \lambda\mathbf{E} [
\widehat{M}_n ] \} \leq C_1 e^{-c_1
\lambda}.
\end{equation}
\end{enumerate}
In particular, if $D$ is a simply connected set containing $0$ and for
all $z \in D$, $\operatorname{dist}(z,D^c) \leq n$, then one can
replace $M_{D',D}$
with $M_D$ in (\ref{eq:emean}), (\ref{eq:emean2}) and (\ref{eq:emean3}).
\end{theorem}
\begin{pf}
The first part follows immediately from Propositions \ref{kmoment} and
\ref{M_Dlb}.

To prove the second part, let $c_0=1/(2C_0)$. Then,
\begin{eqnarray*}
\mathbf{E} \bigl[ \exp\{c_0M_{D',D}/\mathbf{E} [
M_n ] \} \bigr] = \sum
_{k=0}^\infty\frac{(c_0)^k \mathbf{E} [ M_{D',D}^k
]}{k! \mathbf{E} [ M_n ]^k}
\leq\sum_{k=0}^\infty2^{-k} = 2.
\end{eqnarray*}
The third part is then immediate by Markov's inequality.

To prove the last part, we first note that, by Corollary \ref{infdist},
\[
\mathbf{P} \{ \widehat{M}_n > \lambda\mathbf{E} [
\widehat{M}_n ] \} \leq C \mathbf{P} \{ M_{4n}
> \lambda\mathbf{E} [ \widehat{M}_n ] \}.
\]
By Proposition \ref{ExpN} (even though it appears later in this paper,
its proof does not rely on this theorem), $\mathbf{E} [
\widehat{M}_n ] \asymp
n^2 \operatorname{Es}(n)$.
Using Lemma \ref{nEs(n)} and Proposition \ref{M_Dlb}, this implies that
$\mathbf{E} [ \widehat{M}_n ] \asymp\mathbf{E} [ M_{4n}
]$ and, therefore,
\[
\mathbf{P} \{ M_{4n} > \lambda\mathbf{E} [
\widehat{M}_n ] \}
\leq C \mathbf{P} \{ M_{4n} > c \lambda\mathbf{E}
[ M_{4n} ] \}
\leq C e^{-c \cdot c_0 \lambda}
= C_1 e^{-c_1 \lambda}.
\]
\upqed
\end{pf}

\section{Estimating the lower tail of $M_D$ and $\widehat{M}_D$}
\label{lowersec}

\begin{lemma} \label{ProbonLXlow}
There exists $c > 0$ such that the following holds. Suppose that $m$,
$n$, $N$, $K$, $x$, $X$
and $\alpha$ are as in Definition \ref{definition}. Then, for any $z
\in A_n(x)$,
\[
\mathbf{P} \{ z \in\alpha\} \geq c \biggl( \ln
\frac{N}{n} \biggr)^{-3}
\operatorname{Es}(n).
\]
\end{lemma}
\begin{pf}
By Lemma \ref{markovformulalb}, if $Y$ is a random walk started at $x$
conditioned
to hit $z$ before hitting $K$ or leaving $B_N$ and $\tau$ is the last
visit of $z$ before leaving $B_N$, then
\[
\mathbf{P} \{ z \in\alpha\} = G_N^X(x,z) \mathbf
{P} \{ \mathrm{L}(Y[0, \tau]) \cap X^z[1, \sigma_N] =
\varnothing; \mathrm{L}(Y[0, \tau]) \subset B_n(x) \}.
\]
By Lemma \ref{markovgreens}, $G_N^X(x,z) \geq c$. Therefore, if we
imitate the proof
of Proposition \ref{M_Dlb} up to (\ref{eq:M_Dlb4}), it is sufficient
to prove that
for all $v \in\partial B(x;n/16)$, $ | {\arg}(v - x) | \leq
\pi/3$,
%
%
\begin{equation} \label{qw}
\mathbf{P}^{v} \bigl\{ \xi_x < \sigma_{B(x;n/8)} \mid\xi_x < \xi
_K \wedge\sigma_N \bigr\} \geq c \biggl(\ln\frac{N}{n} \biggr)^{-2}
\end{equation}
and for all $w \in\partial B_n(x)$, $ | {\arg}(w - x) |
\leq\pi/3$,
%
%
\begin{equation} \label{df}
\mathbf{P}^{w} \bigl\{ \sigma_N < \xi_{B(x;7n/8)} \mid\sigma_N <
\xi_K \bigr\} \geq c
\biggl(\ln\frac{N}{n} \biggr)^{-1}.
\end{equation}

We first establish (\ref{qw}):
\[
\mathbf{P}^{v} \bigl\{ \xi_x < \sigma_{B(x;n/8)} \mid\xi_x < \xi
_K \wedge\sigma_N \bigr\}
= \frac{\mathbf{P}^{v} \{ \xi_x < \sigma_{B(x;n/8)} \wedge
\xi_K \}}{\mathbf{P}^{v} \{ \xi_x < \xi_K \wedge
\sigma_N \}}.
\]
Let $K' = K \cup\{x\}$. By Lemma \ref{last-exit},
\[
\mathbf{P}^{v} \bigl\{ \xi_x < \sigma_{B(x,n/8)} \wedge\xi_K
\bigr\}
= \frac{G(v,v;B(x;n/8) \setminus K')}{G(v,v; \mathbb{Z}^2\setminus\{
x\})}
\frac{ \mathbf{P}^{x} \{ \xi_v < \xi_{K'} \wedge\sigma
_{B(x;n/8)} \}} {\mathbf{P}^{x} \{ \xi_v < \xi_x
\}}
\]
and
\[
\mathbf{P}^{v} \{ \xi_x < \xi_K \wedge\sigma_N \}
= \frac{G(v,v; B_N \setminus K')}{G(v,v; \mathbb{Z}^2\setminus\{x\})}
\frac{\mathbf{P}^{x} \{ \xi_v < \xi_{K'} \wedge\sigma_{N}
\}}{\mathbf{P}^{x} \{ \xi_v < \xi_x \}}.
\]
Therefore,
\begin{eqnarray*}
&& \mathbf{P}^{v} \bigl\{ \xi_x < \sigma_{B(x;n/8)} \mid\xi_x <
\xi_K \wedge\sigma_N \bigr\} \\
&&\qquad= \frac{G(v,v;B(x;n/8) \setminus K')}{G(v,v; B_N \setminus K')}
\frac
{\mathbf{P}^{x} \{ \xi_v < \xi_{K'} \wedge\sigma_{B(x;n/8)}
\}}{\mathbf{P}^{x} \{ \xi_v < \xi_{K'} \wedge\sigma
_{N} \}}.
\end{eqnarray*}

Since $ | v-x | = n/16$,
\[
G\bigl(v,v; B(x;n/8) \setminus K'\bigr) \geq G\bigl(v,v;B(v; n/16)\bigr
) \geq c \ln n.
\]
Also,
\[
G(v,v; B_N \setminus K') \leq G(v,v; B(v;2N)) \leq C \ln N.
\]
Therefore,
\[
\frac{G(v,v;B(x;n/8) \setminus K')}{G(v,v; B_N \setminus K')} \geq c
\frac{\ln n}{\ln N} \geq c \biggl(\ln\frac{N}{n} \biggr)^{-1}.
\]
To prove (\ref{qw}), it therefore suffices to show that
\[
\mathbf{P}^{x} \{ \xi_v < \xi_{K'} \wedge\sigma_{N}
\} \leq C \ln\frac{N}{n}
\mathbf{P}^{x} \bigl\{ \xi_v < \xi_{K'} \wedge\sigma_{B(x;n/8)}
\bigr\}.
\]
Indeed,
\begin{eqnarray*}
&& \mathbf{P}^{x} \{ \xi_v < \xi_{K'} \wedge\sigma_{N}
\} \\
&&\qquad= \mathbf{P}^{x} \bigl\{ \xi_v < \xi_{K'} \wedge\sigma
_{B(x;n/8)} \bigr\} \\
&&\qquad\quad{} + \sum_{y \in\partial B(x;n/8)} \mathbf{P}^{y} \{ \xi
_v < \xi_{K'} \wedge\sigma_{N} \}\\
&&\qquad\quad\hspace*{57.8pt}{}\times
\mathbf{P}^{x} \bigl\{ S\bigl(\sigma_{B(x;n/8)}\bigr) = y; \sigma_{B(x;n/8)}
< \xi_{K'} \wedge\xi_v \bigr\} \\
&&\qquad\leq\mathbf{P}^{x} \bigl\{ \xi_v < \xi_{K'} \wedge\sigma
_{B(x;n/8)} \bigr\} \\
&&\qquad\quad{} + \sum_{y \in\partial B(x;n/8)} \mathbf{P}^{y} \bigl\{
\xi
_v < \sigma_{B(v;2N)} \bigr\}
\mathbf{P}^{x} \bigl\{ S\bigl(\sigma_{B(x;n/8)}\bigr) = y; \sigma_{B(x;n/8)}
< \xi_{K'} \bigr\}.
\end{eqnarray*}
For all $y \in\partial B(x;n/8)$, $ | y-v | > n/16$ and, thus,
\[
\mathbf{P}^{y} \bigl\{ \xi_v < \sigma_{B(v;2N)} \bigr\} \leq C
\frac{\ln(N/n)}{\ln N}.
\]
Therefore,
\[
\mathbf{P}^{x} \{ \xi_v < \xi_{K'} \wedge\sigma_{N} \}
\leq\mathbf{P}^{x} \bigl\{ \xi_v < \xi_{K'} \wedge\sigma
_{B(x;n/8)} \bigr\}
+ C \frac{\ln(N/n)}{\ln N} \mathbf{P}^{x} \bigl\{ \sigma
_{B(x;n/8)} < \xi_{K'} \bigr\}.
\]
However, by Proposition \ref{dirichlet},
\begin{eqnarray*}
\mathbf{P}^{x} \bigl\{ \sigma_{B(x;n/8)} < \xi_{K'} \bigr\} &\leq&
\mathbf{P}^{x} \bigl\{ \sigma_{B(x;n/16)} < \xi_{K'} \bigr\} \\
&\leq& C \mathbf{P}^{x} \biggl\{ \sigma_{B(x;n/16)} < \xi_{K'} ;
\bigl| {\arg}\bigl(S\bigl(\sigma_{B(x;n/16)}\bigr) - x\bigr) \bigr| \leq
\frac{\pi
}{4} \biggr\} \\
&\leq& C \ln n \mathbf{P}^{x} \bigl\{ \xi_v < \sigma_{B(x;n/8)}
\wedge\xi_{K'} \bigr\}.
\end{eqnarray*}
Thus,
\begin{eqnarray*}
\mathbf{P}^{x} \{ \xi_v < \xi_{K'} \wedge\sigma_{N}
\} &\leq& \biggl(1 + C
\frac{\ln(N/n) \ln n}{\ln N} \biggr)
\mathbf{P}^{x} \bigl\{ \xi_v < \xi_{K'} \wedge\sigma_{B(x;n/8)}
\bigr\} \\
&\leq& C \ln\frac{N}{n} \mathbf{P}^{x} \bigl\{ \xi_v < \xi_{K'}
\wedge\sigma_{B(x;n/8)} \bigr\}.
\end{eqnarray*}

We now prove (\ref{df}):
\[
\mathbf{P}^{w} \bigl\{ \sigma_N < \xi_{B(x;7n/8)} \mid\sigma_N <
\xi_K \bigr\}
= \frac{\mathbf{P}^{w} \{ \sigma_N < \xi_K \wedge\xi
_{B(x;7n/8)} \}}{\mathbf{P}^{w} \{ \sigma_N < \xi_K
\}}.
\]
Let $y_0 \in\partial B_n(x)$ be such that
\[
\mathbf{P}^{y_0} \{ \sigma_N < \xi_K \} = \max_{y \in
\partial B_n(x)} \mathbf{P}^{y} \{ \sigma_N < \xi_K \}.
\]
Then,
\begin{eqnarray*}
\hspace*{-4pt}&& \mathbf{P}^{y_0} \{ \sigma_N < \xi_K \} \\
\hspace*{-4pt}&&\qquad= \mathbf{P}^{y_0} \bigl\{ \sigma_N < \xi_K \wedge
\xi
_{B(x;7n/8)} \bigr\} \\
\hspace*{-4pt}&&\qquad\quad{} + \sum_{u \in\partial_i B(x;7n/8)} \mathbf
{P}^{u} \{
\sigma_N < \xi_K \}
\mathbf{P}^{y_0} \bigl\{ S\bigl(\xi_{B(x;7n/8)}\bigr) = u; \xi
_{B(x;7n/8)} <
\xi_K \wedge\sigma_N \bigr\} \\
\hspace*{-4pt}&&\qquad\leq\mathbf{P}^{y_0} \bigl\{ \sigma_N < \xi_K
\wedge\xi
_{B(x;7n/8)} \bigr\}
+ \mathbf{P}^{y_0} \{ \sigma_N < \xi_K \} \mathbf
{P}^{y_0} \bigl\{ \xi_{B(x;7n/8)} < \xi_K \wedge\sigma_N \bigr\}
\\
\hspace*{-4pt}&&\qquad\leq\mathbf{P}^{y_0} \bigl\{ \sigma_N < \xi_K
\wedge\xi
_{B(x;7n/8)} \bigr\}
+ \mathbf{P}^{y_0} \{ \sigma_N < \xi_K \} \mathbf
{P}^{y_0} \bigl\{ \xi_{B(x;7n/8)} < \sigma_{B(x;2N)} \bigr\}.
\end{eqnarray*}
However, by Proposition \ref{p:potential},
\[
\mathbf{P}^{y_0} \bigl\{ \xi_{B(x;7n/8)} < \sigma_{B(x;2N)}
\bigr\} \leq1 - \frac{c}{\ln(N/n)}
\]
and, therefore,
\[
\mathbf{P}^{y_0} \{ \sigma_N < \xi_K \} \leq C \ln
\frac{N}{n} \mathbf{P}^{y_0} \bigl\{ \sigma_N < \xi_K \wedge\xi
_{B(x;7n/8)} \bigr\}.
\]
This establishes (\ref{df}) for the special case where $w = y_0$.
However, we can apply Lemma \ref{reflection} twice, as in Corollary
\ref{reflectcor}, to conclude that
\[
\mathbf{P}^{w} \bigl\{ \sigma_N < \xi_K \wedge\xi_{B(x;7n/8)}
\bigr\} \geq c \max_{y
\in\partial B_n(x)}
\mathbf{P}^{y} \bigl\{ \sigma_N < \xi_K \wedge\xi_{B(x;7n/8)}
\bigr\}.
\]
Therefore,
\[
\frac{\mathbf{P}^{w} \{ \sigma_N < \xi_K \wedge\xi
_{B(x;7n/8)} \}}{\mathbf{P}^{w} \{ \sigma_N < \xi_K
\}}
\geq c \frac{\mathbf{P}^{y_0} \{ \sigma_N < \xi_K \wedge\xi
_{B(x;7n/8)} \}}{\mathbf{P}^{y_0} \{ \sigma_N < \xi_K
\}}
\geq c \biggl( \ln\frac{N}{n} \biggr)^{-1}.
\]
\upqed
\end{pf}
\begin{prop} \label{ExpN}
\begin{enumerate}
\item There exists $C < \infty$ such that for any $m$, $n$, $N$, $K$
and $x$ as in Definition \ref{definition},
\[
C^{-1} \biggl(\ln\frac{N}{n} \biggr)^{-3} n^2 \operatorname{Es}(n)
\leq\mathbf{E} [ M_{m,n,N,x}^K ] \leq C \biggl(\ln
\frac{N}{n} \biggr)
n^2 \operatorname{Es}(n) .
\]
\item
\[
\mathbf{E} [ M_n ] \asymp\mathbf{E} [
\widehat{M}_n ] \asymp n^2 \operatorname{Es}(n).
\]
\end{enumerate}
\end{prop}
\begin{pf}
We first prove part 1. Let $\alpha$ be as in Definition \ref
{definition}. Then, by Lem\-ma~\ref{ProbonLXlow},
\begin{eqnarray*}
\mathbf{E} [ M_{m,n,N,x}^K ] &=&
\sum_{z \in A_n(x)} \mathbf{P} \{ z \in\alpha\}
\geq\sum_{z \in A_n(x)}
c \biggl(\ln\frac{N}{n} \biggr)^{-3} \operatorname{Es}(n) \\
&\geq& c
\biggl(\ln\frac
{N}{n} \biggr)^{-3} n^2 \operatorname{Es}(n).
\end{eqnarray*}
To prove the other direction, note that by Proposition \ref{kpoints},
with $k=1$, for any $z \in\alpha$,
\[
\mathbf{P} \{ z \in\alpha\} \leq
\mathbf{P} \{ z \in\mathrm{L} (X[0,\sigma_N])
\} \leq C G^X_N(x,z) \operatorname{Es}(n).
\]
By Lemma \ref{markovgreens}, $G^X_N(x,z) \leq C \ln(N/n)$ and, therefore,
\[
\mathbf{E} [ M_{m,n,N,x}^K ] =
\sum_{z \in A_n(x)} \mathbf{P} \{ z \in\alpha\}
\leq C \biggl(\ln\frac
{N}{n} \biggr) n^2 \operatorname{Es}(n).
\]

We now prove part 2. The fact that $\mathbf{E} [ M_n ]
\asymp n^2 \operatorname{Es}(n)$
follows immediately from Theorem \ref{kmoment} and Proposition \ref{M_Dlb}.

In order to show that $\mathbf{E} [ \widehat{M}_n ]
\asymp n^2 \operatorname{Es}(n)$, let
$\beta$ be $\mathrm{L}(S[0,\sigma_{4n}])$ from $0$ up to its
first exit from
the ball $B_n$. By Corollary \ref{infdist}, $\beta$ has the same
distribution, up to constants, as $\widehat{S}[0,\widehat{\sigma
}_n]$ and thus it
suffices to show that
\[
\sum_{z \in B_n} \mathbf{P} \{ z \in\beta\} \asymp
n^2 \operatorname{Es}(n).
\]

To begin with,
\[
\sum_{z \in B_n} \mathbf{P} \{ z \in\beta\}
\leq\sum_{z \in B_{4n}} \mathbf{P} \{ z \in\mathrm{L}(S[0,\sigma_{4n}])
\} \asymp
n^2 \operatorname{Es}(4n).
\]
By Lemma \ref{Esdecreasing}, the latter is less than a constant times
$n^2 \operatorname{Es}(n)$.

To prove the other direction, the number of steps of $\beta$ is
strictly larger
than $M_{n,m,N,x}^K$, where $m=0$, $N=4n$, $x=0$ and $K = \varnothing$.
Therefore, by part 1 and Lemma \ref{Es(kn)},
we have
\[
\sum_{z \in B_n} \mathbf{P} \{ z \in\beta\} \geq
\mathbf{E} [ M_{n,0,4n,0}^\varnothing] \geq c n^2
\operatorname{Es}(4n)
\geq c n^2 \operatorname{Es}(n).
\]
\upqed
\end{pf}
\begin{prop} \label{ExpN^2} There exists $C < \infty$ such that if
$m$, $n$, $N$, $K$ and $x$ are as in Definition \ref{definition},
then
\[
\mathbf{E} [ (M_{m,n,N,x}^K)^2 ] \leq C \biggl(\ln
\frac{N}{n} \biggr)^2
n^4 \operatorname{Es}(n)^2.
\]
\end{prop}
\begin{pf}
Let $\alpha$ be as in Definition \ref{definition}. Then, by
Proposition \ref{kpoints},
\begin{eqnarray*}
\mathbf{E} [ (M_{m,n,N,x}^K)^2 ] &=& \mathbf{E} \biggl[ \biggl(
\sum_{z \in A_n(x)} \mathbh{1} _{\{z \in\alpha\}}
\biggr)^2 \biggr] \\
&=& \sum_{z,w \in A_n(x)} \mathbf{P} ( z,w \in\alpha
) \\
&\leq& C \sum_{z,w \in A_n(x)} G^X_N(x,z)G^X_N(z,w) \operatorname
{Es}(r_z) \operatorname{Es}(r_w),
\end{eqnarray*}
where $r_z = \operatorname{dist}(z,\partial B_N) \wedge| z-x
| \wedge| z-w |$ and
$r_w = \operatorname{dist}(z,\partial B_N) \wedge| z-w
|$. However, since $z$
and $w$ are in $A_n(x)$, $r_z$
and $r_w$ are comparable to $ | z-w |$. Therefore, by
Lemmas \ref{markovgreens}, \ref{Es(kn)} and the fact that
\[
G^X_N(z,w) = G_{B_N \setminus K}(z,w) \frac{\mathbf{P}^{w} \{
\sigma_N < \xi_K \}}
{\mathbf{P}^{z} \{ \sigma_N < \xi_K \}} \leq C
G_{B_{2N}(z)}(z,w) \leq C \ln
\frac{2N}{ | z-w |},
\]
we have
\begin{eqnarray*}
\mathbf{E} [ (M_{m,n,N,x}^K)^2 ] &\leq& C \ln\frac
{N}{n} \sum_{z,w \in
A_n(x)} \ln\frac{2N}{ | z-w |} \operatorname{Es}( |
z-w |)^2 \\
&\leq& C \ln\frac{N}{n} \sum_{z \in A_n(x)} \sum_{w \in B_n(z)} \ln
\frac{2N}{ | z-w |} \operatorname{Es}( | z-w
|)^2 \\
&\leq& C \ln\frac{N}{n} \sum_{z \in A_n(x)} \sum_{k = 1}^n k \ln
\frac{N}{k} \operatorname{Es}(k)^2 \\
&\leq& C \ln\frac{N}{n} n^2 \Biggl(\sum_{k = 1}^n k \ln\frac{n}{k}
\operatorname{Es}(k)^2 + \sum_{k = 1}^n k \ln\frac{N}{n}
\operatorname{Es}(k)^2 \Biggr).
\end{eqnarray*}
By Corollary \ref{Es_sums}, both of the sums above are bounded by $C
\ln(N/n) n^2 \operatorname{Es}(n)^2$ which completes the proof.
\end{pf}
\begin{cor} \label{secondcor}
There exist $C < \infty$ and $c_2, c_3 > 0$ such that if $m$, $n$,
$N$, $K$ and $x$ are as in Definition \ref{definition},
then:
\begin{enumerate}
\item
\[
\mathbf{E} [ (M_{m,n,N,x}^K)^2 ] \leq C \biggl( \ln
\frac{N}{n} \biggr)^{8}
\mathbf{E} [ (M_{m,n,N,x}^K) ]^2;
\]
\item
\[
\mathbf{P} \biggl\{ M_{m,n,N,x}^K \leq c_2 \biggl(\ln\frac
{N}{n} \biggr)^{-3} \mathbf{E} [ M_n ] \biggr\} \leq
1 - c_3 \biggl(\ln\frac{N}{n} \biggr)^{-8}.
\]
\end{enumerate}
\end{cor}
\begin{pf}
The first part follows immediately from Propositions \ref{ExpN} and
\ref{ExpN^2}.

To prove the second part, by a standard second moment result (see,
e.g., \cite{LL08}, Lemma 12.6.1), for any $0 < r < 1$,
\[
\mathbf{P} \{ M_{m,n,N,x}^K \leq r \mathbf{E} [
M_{m,n,N,x}^K ] \}
\leq1 - \frac{(1-r)^2\mathbf{E} [ M_{m,n,N,x}^K
]^2}{\mathbf{E} [ (M_{m,n,N,x}^K)^2 ]}.
\]
Letting $r=1/2$ and using part 1, one obtains that
\[
\mathbf{P} \biggl\{ M_{m,n,N,x}^K \leq\frac{1}{2} \mathbf{E} [
M_{m,n,N,x}^K ] \biggr\} \leq1
- c_3 \biggl(\ln\frac{N}{n} \biggr)^{-8}.
\]
Finally, by again using Proposition \ref{ExpN}, we have
\[
\mathbf{E} [ M_{m,n,N,x}^K ] \geq c \biggl(\ln\frac
{N}{n} \biggr)^{-3} \mathbf{E} [ M_n ].
\]
\upqed\end{pf}
\begin{lemma} \label{expectationrelation}
For all $\varepsilon> 0$, there exist $C(\varepsilon) < \infty$ and
$N(\varepsilon) < \infty$ such that for
all $n \geq N(\varepsilon)$ and $k \geq1$,
\[
\mathbf{E} [ M_{kn} ] \leq C(\varepsilon)
k^{5/4+\varepsilon} \mathbf{E} [ M_n ]
\]
and
\[
\mathbf{E} [ \widehat{M}_{kn} ] \leq C(\varepsilon)
k^{5/4+\varepsilon} \mathbf{E} [ \widehat{M}_n ].
\]
\end{lemma}
\begin{Remark*}
It is possible to take
$\varepsilon= 0$ in the inequality above, but, in that case, $N$ has to
depend on $k$.
\end{Remark*}
\begin{pf*}{Proof of Lemma \ref{expectationrelation}}
The second statement follows immediately from the first, by
Proposition \ref{ExpN}.

By Proposition \ref{ExpN} and Theorem \ref{bigthm}, part 3, we have
\[
\mathbf{E} [ M_{kn} ] \leq C (kn)^2 \operatorname
{Es}(kn) \leq C (kn)^2 \operatorname{Es}(n) \operatorname{Es}(n,kn).
\]
By Lemma \ref{agrt}, there exist $C(\varepsilon) < \infty$ and
$N(\varepsilon)$ such that for all $n \geq N(\varepsilon)$,
\[
\operatorname{Es}(n,kn) \leq C(\varepsilon) k^{-3/4+\varepsilon}.
\]
Therefore,
\[
\mathbf{E} [ M_{kn} ] \leq C(\varepsilon)
k^{5/4+\varepsilon} n^2 \operatorname{Es}(n).
\]
Finally, by a second application of Proposition \ref{ExpN},
we obtain
\[
n^2 \operatorname{Es}(n) \leq C \mathbf{E} [ M_n ].
\]
\upqed\end{pf*}
\begin{prop}[(See Figure \ref{iteratepic})] \label{iterate}
\begin{enumerate}
\item Let $c_2$ be as in Corollary \ref{secondcor}. There then exists
$c_4 > 0$ such that for all $n$ and all $k \geq2$,
\[
\mathbf{P} \{ M_{kn} \leq c_2 ( \ln k )^{-3}
\mathbf{E} [ M_n ] \} \leq
e^{-c_4 k (\ln k)^{-8}}.
\]
\item There exist $c_5, c_6 > 0$ and $C < \infty$ such that for all
$n$ and $k \geq2$,
\[
\mathbf{P} \{ \widehat{M}_{kn} \leq c_5 ( \ln k
)^{-3} \mathbf{E} [ \widehat{M}_n ] \} \leq C
e^{-c_6 k (\ln k)^{-8}}.
\]
\end{enumerate}
\end{prop}
\begin{pf}
We first prove part 1.

%
%
\begin{figure}

\includegraphics{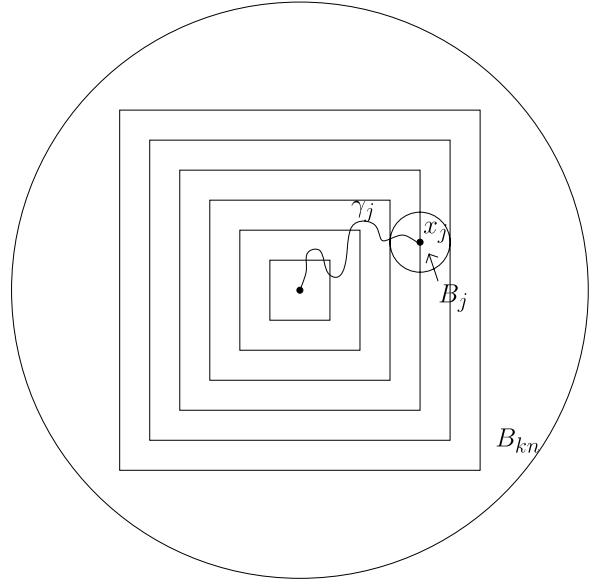}

\caption{The setup for Proposition \protect\ref{iterate}.}
\label{iteratepic}
\end{figure}

Let $k'= \lfloor k/\sqrt{2} \rfloor$. Then, $R_{k'n} \subset B_{kn}$.
We view the loop-erased random walk $\mathrm{L}(S[0,\sigma
_{kn}])$ as a
distribution on the
set $\Omega_{kn}$ of self-avoiding paths $\gamma$ from the origin to
$\partial B_{kn}$.
Given such a $\gamma$, let $\gamma_j$ be its restriction from $0$ to
the first exit of
$R_{jn}$, $j=0, \ldots, k'$. Let $\mathcal{F}_j$ be the $\sigma
$-algebra generated by
the $\gamma_j$. For $j=0, \ldots, k'-1$, let $x_j(\gamma) \in
\partial R_{jn}$ be the point
where $\gamma$ first exits $R_{jn}$ and $B_j = B_{n}(x_j)$. Finally,
for $j=1, \ldots, k'$,
let $\alpha_{j}(\gamma)$ be $\gamma$ from $x_{j-1}$ up to the first
exit of $B_{j-1}$ and
let $N_j(\gamma)$ be the number of steps of $\alpha_j$ in
$A_n(x_{j-1})$ [where $A_n(x)$ is as in Definition \ref{definition}].
Note that $N_j \in\mathcal{F}_{j}$.

Then,
\begin{eqnarray*}
&& \mathbf{P} \{ M_{kn} \leq c_2 ( \ln k )^{-3}
\mathbf{E} [ M_n ] \} \\
&&\qquad\leq\mathbf{P} \Biggl\{ \sum_{j=1}^{k'} N_j \leq c_2 ( \ln
k )^{-3} \mathbf{E} [ M_n ] \Biggr\} \\
&&\qquad\leq\mathbf{P} \Biggl( \bigcap_{j=1}^{k'} \{ N_j \leq c_2
( \ln k )^{-3} \mathbf{E} [ M_n ]
\} \Biggr) \\
&&\qquad= \mathbf{E} \Biggl[ \Biggl( \prod_{j=1}^{k'-1} \mathbh{1}_{\{
N_j \leq c_2 ( \ln k )^{-3} \mathbf{E} [ M_n
] \}} \Biggr) \mathbf{P} \{ N_{k'} \leq c_2 ( \ln
k )^{-3} \mathbf{E} [ M_n ] \mid\mathcal
{F}_{k'-1} \} \Biggr].
\end{eqnarray*}

However, by the domain Markov property, for all $j = 1, \ldots, k'$,
\[
\mathbf{P} \{ N_j \leq c_2 ( \ln k )^{-3} \mathbf
{E} [ M_n ] \mid\mathcal{F}_{j-1} \}(\gamma)
= \mathbf{P} \bigl\{ M_{jn, n, kn, x_j(\gamma)}^{\gamma_j} \leq
c_2 ( \ln k )^{-3} \mathbf{E} [ M_n ]
\bigr\}.
\]
Furthermore, by Corollary \ref{secondcor},
\[
\mathbf{P} \bigl\{ M_{jn, n, kn, x_j(\gamma)}^{\gamma_j} \leq
c_2 ( \ln k )^{-3} \mathbf{E} [ M_n ]
\bigr\}
\leq1 - c_3(\ln k)^{-8}.
\]
Therefore, by applying the above inequality $k'$ times, we obtain
\begin{eqnarray*}
\mathbf{P} \{ M_{kn} \leq c_2 ( \ln k )^{-3}
\mathbf{E} [ M_n ] \} &\leq\bigl(1 - c_3(\ln
k)^{-8} \bigr)^{k'}
\leq e^{-c_3(\ln k)^{-8}k'}.
\end{eqnarray*}

The proof of part 2 is analogous.
By Proposition \ref{ExpN}, it suffices to show that
\[
\mathbf{P} \{ \widehat{M}_{kn} \leq c_2 ( \ln k
)^{-3} \mathbf{E} [ M_n ] \} \leq e^{-c_6k (\ln
k)^{-8}}.
\]
However, by Corollary \ref{infdist}, $\widehat{S}[0,\widehat{\sigma}_{kn}]$
has the same
distribution, up to constants, as $\mathrm{L}(S[0,\sigma
_{4kn}])$ from $0$
up to its first exit
of the ball $B_{kn}$. Therefore, we can apply the previous iteration
argument to obtain that
\[
\mathbf{P} \{ \widehat{M}_{kn} \leq c_2 ( \ln k
)^{-3} \mathbf{E} [ M_n ] \}
\leq C \bigl(1 - c_3(\ln4k)^{-8} \bigr)^{k'} \leq C e^{-c_6k (\ln k)^{-8}}.
\]
\upqed\end{pf}
\begin{theorem} \label{main}
For all $\varepsilon> 0$, there exist
$C_2(\varepsilon) <\infty$, $C_3(\varepsilon) < \infty$,
$c_7(\varepsilon
)>0$ and
$c_8(\varepsilon) > 0$ such that for all $\lambda> 0$ and all $n$:
\begin{enumerate}
\item
\[
\mathbf{P} \{ \widehat{M}_n < \lambda^{-1} \mathbf{E} [ \widehat
{M}_n ] \}
\leq C_2(\varepsilon) e^{-c_7(\varepsilon)\lambda^{4/5 - \varepsilon}};
\]
\item for all $D \supset B_n$, $\lambda> 0$,
\[
\mathbf{P} \{ M_D < \lambda^{-1} \mathbf{E} [ M_n
] \}
\leq C_3(\varepsilon) e^{-c_8(\varepsilon) \lambda^{4/5 -
\varepsilon}}.
\]
\end{enumerate}
\end{theorem}
\begin{pf}
The second part follows from the first since, by Corollary
\ref{infdist}, Proposition \ref{ExpN}
and Lemma \ref{expectationrelation}, we have
\begin{eqnarray*}
\mathbf{P} \{ M_D < \lambda^{-1} \mathbf{E} [ M_n
] \} &\leq& C \mathbf{P} \{ \widehat{M}_{n/4} <
\lambda^{-1} \mathbf{E} [ M_n ] \} \\
&\leq& C \mathbf{P} \{ \widehat{M}_{n/4} < C \lambda^{-1}
\mathbf{E} [ \widehat{M}_n ] \} \\
&\leq& C \mathbf{P} \{ \widehat{M}_{n/4} < C \lambda^{-1}
\mathbf{E} [ \widehat{M}_{n/4} ] \}.
\end{eqnarray*}

We now prove the first part. We will prove the result for all
$\varepsilon
$ such
that $0 < \varepsilon< 7/40$ and note that for such $\varepsilon$,
\[
\frac{5}{4} + \varepsilon\leq\frac{1}{4/5 - \varepsilon} \leq
\frac
{5}{4} + 2 \varepsilon.
\]
Clearly, this will imply that the result holds for all $\varepsilon> 0$.

Fix such an $\varepsilon> 0$. We will show that there exist $C <
\infty
$, $c_7 > 0$, $\lambda_0$
and $N$ such that, for $\lambda> \lambda_0$ and $n \geq N$,
%
%
\begin{equation} \label{eq:main2}
\mathbf{P} \{ \widehat{M}_n < \lambda^{-1} \mathbf{E} [ \widehat
{M}_n ] \} \leq C
e^{-c_7\lambda^{4/5 - \varepsilon}}.
\end{equation}
We claim that this implies the statement of the theorem with
\[
C_2 = C \vee e^{4c_7 (\lambda_0 \vee N)^{4/5 - \varepsilon}}.
\]
To see this, if $\lambda< \lambda_0$, then, for any $n$,
\[
\mathbf{P} \{ \widehat{M}_n < \lambda^{-1} \mathbf{E} [ \widehat
{M}_n ] \} \leq1 \leq C_2
e^{-c_7 \lambda^{4/5 - \varepsilon}}.
\]
Next, if $n \leq N$, then,
for any $\lambda$,
\[
\mathbf{P} \{ \widehat{M}_n < \lambda^{-1} \mathbf{E} [ \widehat
{M}_n ] \} \leq\mathbf{P} \{
\widehat{M}_n < 4 \lambda^{-1} n^2 \}
\]
since $\mathbf{E} [ \widehat{M}_n ] \le|B_n| <
4n^2$. If $\lambda> 4 n$, then
the above probability is $0$
since $\mathbf{P} \{ \widehat{M}_n \geq n \} = 1$. If
$\lambda< 4 n \leq4 N$, then
\[
C_2 e^{-c_7 \lambda^{4/5 - \varepsilon}} \geq e^{4 c_7 N^{4/5 -
\varepsilon
}} e^{-4 c_7 N^{4/5 - \varepsilon}} = 1.
\]

We now prove (\ref{eq:main2}). Let $c_5$ be as in Proposition
\ref{iterate}, and
$C^* = C(\varepsilon/2)$ and $N_0=N(\varepsilon/2)$ be as in Lemma
\ref{expectationrelation}.
Let
\[
k = c_5 (C^*)^{-1} \lambda^{4/5 - \varepsilon/2 }.
\]
We choose $\lambda_0$ so that for all $\lambda> \lambda_0$, $k \geq2$,
$k^{\varepsilon/2} > (\ln k)^3$ and $k (\ln k)^{-8} \ge\lambda
^{4/5-\varepsilon}$.
We also choose $N = 4 N_0^{5}$. Then, for all $n \geq N$ and $\lambda>
\lambda_0$,
%
%
\begin{equation} \label{eq:main}
\mathbf{E} [ \widehat{M}_{kn} ] \leq C^*
k^{5/4+\varepsilon/2} \mathbf{E} [ \widehat{M}_n ]
\leq c_5 k^{-\varepsilon/2} \lambda\mathbf{E} [ \widehat
{M}_n ].
\end{equation}

First, suppose that $n/k \leq N_0$. Then,
\[
\lambda^{-1} \leq k^{-5/4} \leq(N_0n^{-1})^{5/4} \leq1/(4n)
\]
and so $\lambda^{-1} \mathbf{E} [ \widehat{M}_n ]
\le n$. Hence, since $\widehat
{M}_n \geq n$ almost surely,
\[
\mathbf{P} \{ \widehat{M}_n < \lambda^{-1} \mathbf{E} [ \widehat
{M}_n ] \}
\leq\mathbf{P} \{ \widehat{M}_n < n \} = 0.
\]

If $n/k \geq N_0$, then, by (\ref{eq:main}) and Proposition \ref{iterate},
\begin{eqnarray*}
\mathbf{P} \{ \widehat{M}_n < \lambda^{-1} \mathbf{E} [ \widehat
{M}_n ] \}
&=& \mathbf{P} \bigl\{ \widehat{M}_{k(n/k)} < \lambda^{-1}
\mathbf{E} \bigl[ \widehat{M}_{k(n/k)} \bigr] \bigr\}
\\
&\leq&\mathbf{P} \bigl\{ \widehat{M}_{k(n/k)} < c_5
k^{-\varepsilon/2} \mathbf{E} [ \widehat{M}_{n/k} ]
\bigr\} \\
&\leq&\mathbf{P} \bigl\{ \widehat{M}_{k(n/k)} < c_5 (\ln k)^{-3}
\mathbf{E} [ \widehat{M}_{n/k} ] \bigr\} \\
&\leq& C e^{-c_6k(\ln k)^{-8}} \\
&\leq& C e^{-c_7\lambda^{4/5 - \varepsilon}}.
\end{eqnarray*}
\upqed
\end{pf}

%

%
\printaddresses

\end{document}